\documentclass[1 leqno,11pt,psfig]{amsart}
\usepackage{amssymb, amsmath,amsmath,latexsym,amssymb,amsfonts,amsbsy, amsthm,mathtools,graphicx,CJKutf8,CJKnumb,CJKulem,color}
\usepackage{graphicx,epsfig}
\usepackage{graphicx}

\usepackage{amsmath}
\usepackage{mathrsfs}
\usepackage{amssymb}
\usepackage{cite}

\usepackage{amssymb,mathrsfs,amsfonts,amsmath,amsbsy,tikz}\usepackage{fontenc}\usepackage{textcomp}

\numberwithin{equation}{section}

\allowdisplaybreaks

\def\C{\mathbb{C}}

\def\Ls{\mathscr{L}}
\def\A{\mathbf{A}}

\def\I{\mathcal{I}}

\let\ld=\lambda

\let\De=\Delta

\let\af=\alpha
\let\bt=\beta

\newcommand{\ba}{\begin{array}}
\newcommand{\ea}{\end{array}}
\newcommand{\be}{\begin{equation}}
\newcommand{\ee}{\end{equation}}
\newcommand{\ban}{\begin{eqnarray*}}
\newcommand{\ean}{\end{eqnarray*}}

\numberwithin{equation}{section}

\newtheorem{theorem}{Theorem}[section]
\newtheorem{definition}[theorem]{Definition}
\newtheorem{lemma}[theorem]{Lemma}

\newtheorem{corollary}[theorem]{Corollary}
\newtheorem{remark}[theorem]{Remark}

\begin{document}

\title[Quasi self-adjoint extensions]
{Quasi self-adjoint extensions of  a class of formally non-self-adjoint discrete Hamiltonian systems}

\author{Guojing Ren}
\address{School of Mathematics and Quantitative Economics,
 Shandong  University of Finance and Economics, Jinan, Shandong 250014, P. R.
China}
\email{gjren@sdufe.edu.cn}

\author{Guixin Xu}
\address{School of Mathematics and Statistics, Beijing Technology and Business University, Beijing, P. R. China}
\email{guixinxu$_-$ds@163.com}

\maketitle

\begin{abstract}
This paper is concerned with the characterizations of quasi self-adjoint extensions of
a class of formally non-self-adjoint discrete Hamiltonian systems.
Some properties of the solutions and the characterization of the
minimal linear relations of the non-self-adjoint systems  are obtained.
A bijective projection between all the quasi self-adjoint extensions of  non-self-adjoint  systems and
all the self-adjoint extensions of the self-adjoint  systems  generated by the  non-self-adjoint  Hamiltonian systems
is established in the general case.
When the system is in the limit point case and $\I=[a,\infty)$,
a  complete characterization of  all the quasi self-adjoint extensions  is obtained
by a  subspace $Q\subset \C^{2n}$ with $\dim Q=n$ in terms of boundary conditions.
\end{abstract}

{\bf \it Keywords}:\ discrete Hamiltonian system, non-self-adjoint, quasi self-adjoint extension, deficiency index, limit point case.

{2020 {\bf \it Mathematics Subject Classification}}:   39A70, 47B28

\section{Introduction}

In this paper we consider the  following  discrete linear Hamiltonian system
\begin{align}\label{HS1}
    (\Ls_1 y)(t):=J\Delta y(t)-P(t)R(y)(t)=\ld W(t))R(y)(t),\quad t \in \I,
\end{align}
where  $\I:=\{t\}_{t=a}^{b}$ is an integer interval,   $a$ is a finite integer or $a=-\infty$,
and $b$ is a finite integer or $b=+\infty$, $b-a\ge 1$;   $J$ is the  canonical symplectic matrix, i.e.,
\begin{eqnarray*}
J=\left(\begin{array}{cc} 0&-I_n\\I_n&0\end{array}\right),
\end{eqnarray*}
and $I_n$ is the $n\times n$  identity matrix; $\De$ is the forward difference operator, i.e.,
$\De y(t)=y(t+1)-y(t)$; $W(t)$ and $P(t)$
are $2n\times 2n$ complex-valued matrices defined on $\I$, and the weight function $W(t)\ge 0$;
the partial right shift operator $R(y)(t)=(u^T(t+1), v^T(t))^T$ with $y(t)=(u^T(t), v^T(t))^T$ and
$u(t)$, $v(t)\in \C^n$;  $\ld\in \C$ is a  spectral parameter.

System (1.1) can be regarded as a discretization of the following continuous Hamiltonian system:
\begin{align}\label{conti HS}
    Jy'(t)=(P(t)+\ld W(t))y(t),\quad t \in (a,b).
\end{align}
When $P$ is Hermitian, that is $P^*=P$,  system (\ref{conti HS}) and system (1.1)
are  said to be formally self-adjoint continuous and discrete Hamiltonian systems, respectively.
Otherwise, they are said to be formally non-self-adjoint.
Here we point out that although many of their properties of (1.1) and (\ref{conti HS}) are similar,
some methods used in the discrete case are different from those used in the continuous case,
and it seems more complicated and difficult in the discrete case.
See \cite{Atkinson1964,Clark2002,Dunford1963,Lesch2003,Naimark1968,Sun2010,Wang2009,Weidmann1980}
for the continuous case and see \cite{Ren2011,Ren2014AML,Ren2014Laa,Shi2006} for the discrete case.

Characterizations of self-adjoint extensions are most fundamental in
the study of spectral problems for  formally self-adjoint differential and difference equations,
and they have been widely discussed.
The corresponding maximal and minimal operators play a key role in the study of  self-adjoint extensions.
It has been shown in \cite{Ren2014AML} that the maximal operator generated by (1.1) with $P^*=P$
is multi-valued in any non-trivial case,
and the minimal operator  may be multi-valued or not densely defined in general.
Therefore, the classical theory of linear operators can not be applied to (1.1).
To deal with these problems, the theory of linear relations has been  developed.
The readers are refereed to \cite{Arens1961,Coddington1973,Cross1998, Lesch2003, Shi2012, XuRen2024} for details.

We recall some important existing results for (1.1) with $P^*=P$.
By  $H$, $H_{00}$, and $H_{0}$ denote the maximal, pre-minimal,
and minimal linear relations corresponding to (1.1) with $P^*=P$.
It has been shown by \cite[Lemma 3.1]{Ren2011} that $H_0^*=H$.
By $d_{\pm}(H_0)$ denote the  positive and negative deficiency  indices  of  $H_0$.
When $\I=[a,\infty)$, a classification of Hamiltonian systems is given in terms of the values of $d_{\pm}(H_0)$.
System (1.1) is said to be in the limit point case if $d_{\pm}(H_0)=n$ and to be in the limit circle case if $d_{\pm}(H_0)=2n$
\cite[Section 5]{Shi2006}.
With the additional assumption that $d_{+}(H_0)=d_{-}(H_0)$,
a complete characterization of all the self-adjoint  extensions for system (1.1) are obtained in
terms of boundary conditions via linear independent square summable
solutions \cite[Section 5]{Ren2014Laa}.

However, there are an increasing number of problems in physics which require the analysis
of spectral problems of formally non-self-adjoint differential and difference equations.
For example, second order linear differential and difference equations with complex valued potential term, and
more generally complex symmetric operators have many applications \cite{Garcia2014}.
Many researchers have discussed these problems and some good results have been obtained in recent years.
See \cite{Brown2003,Brown2019,Cascaval2004,Locker2000} for the continuous case
and \cite{Monaquel2007,Sun2011AMC,Wilson2005,Zhang2023} for the discrete case.

Our  goal of this manuscript is to  establish the characterizations of the quasi selfadjoint
extension of the minimal linear relation generated by (1.1).
To do this, we consider the adjoint system of (1.1):
\begin{align}\label{HS2}
    (\Ls_2 y)(t):=J\Delta y(t)-P^*(t)R(y)(t)=\ld W(t))R(y)(t),\quad t \in \I.
\end{align}

Throughout the whole paper, we assume that the weight
function $W(t)$ is of the block diagonal form,
\begin{eqnarray*}
W(t)={\rm diag }\{W_1(t),W_2(t)\},
\end{eqnarray*}
where $W_j(t)\ge 0$ are $n\times n$ Hermitian matrices, $j=1,2$.
Let  $P(t)$ be blocked as
\begin{eqnarray}\label{P=A,B,C,D}
P(t)=\begin{pmatrix}-C(t)& D(t)\\A(t) &B(t)\end{pmatrix},
\end{eqnarray}
where $A(t)$, $B(t)$, $C(t)$ and $D(t)$ are $n\times n$ complex-valued matrices.
Then  $(1.1)$ can be written as
\begin{eqnarray}\label{HS blocked}
\left\{\begin{array}{l}
         \De u(t)=A(t)u(t+1)+(B(t)+\ld W_2(t))v(t),\\
         \De v(t)=(C(t)-\ld W_1(t))u(t+1)-D(t)v(t),\quad t\in \I.
       \end{array}\right.
\end{eqnarray}

The  main assumptions in this manuscript are:

\begin{itemize}
\item[${\bf(A_1)}$] $I_n-A(t)$ and $I_n-D(t)$ are invertible on $\I$.
\item [\rm {($\bf A_{2}$)}]  There exists a finite subinterval
$\I_0\subseteq \I$ such that for any $g_i$ and $y_i$ with $y_i(k_i)\neq 0$ for some $k_i\in \I$,
which satisfy $(\Ls_iy_i)(t)=W(t)R(g_i)(t)$ on $\I$,  $i=1,2$,
it follows that
\begin{align*}
\sum_{t\in \I_0}R(y_i)^*(t)W(t)R(y_i)(t)>0.
\end{align*}
\end{itemize}

Assumption $(A_1)$ is required to ensure the existence and uniqueness of the solution of any
initial value problem for $(1.1)$. For $i=1,2$, by  $H_{i}$ and $H_{i,0}$ denote the maximal
and the minimal linear relations corresponding to (1.1) and (\ref{HS2}), respectively.
By some technique we combine (1.1) and (\ref{HS2}) into  a formally self-adjoint Hamiltonian system $(\ref{HS3})$.
Based on the results on formally self-adjoint Hamiltonian systems, we establish a bijective projection between all
the quasi self-adjoint extensions of $\{H_{1,0},H_{2,0}\}$ and
all the self-adjoint extensions of the minimal linear relations $\mathbf{H}_0$  generated by $(\ref{HS3})$.

\smallskip

{\bf Theorem\;  \ref{q self 1}.}
{\it Assume that $(A_1)$ and $(A_2)$ hold.
\begin{enumerate}
  \item If $T$ is a quasi self-adjoint extension of  $\{H_{1,0},H_{2,0}\}$,
then  $T^*$ is a quasi self-adjoint extension of  $\{H_{2,0},H_{1,0}\}$,
and  $\mathbf{T}$ generated by $\{T, T^*\}$ is a  self-adjoint extension of $\mathbf{H}_0$.
\item If $\mathbf{T}$  is a  self-adjoint extension of $\mathbf{H}_0$, t
hen $\mathbf{T}$ is generated by $\{T, T^*\}$,
where $T$ is a  quasi self-adjoint extension of  $\{H_{1,0},H_{2,0}\}$.
\end{enumerate}}
\smallskip

In the case that $\I=[a,\infty)$ and $(\ref{HS4})$ is in limit point case,
a complete characterization of all the extensions and all the quasi self-adjoint  extensions for the system (1.1) is obtained in
terms of boundary conditions.
\smallskip

{\bf Theorem\; \ref{Q and H_1}.}
{\it Let $\I=[a,\infty)$. Assume  that  $(A_1)$ and $(A_{2})$ hold,  and (\ref{HS4}) is in limit point case. Then
\begin{enumerate}
  \item $T\in Ext\{H_{1,0},H_{2,0}\}$ if and only if there exists $Q\subset \C^{2n}$ such that
\begin{equation*}
     T= \{(y,[g])\in H_{1}:\; y(a)\in Q\}
  \end{equation*}
  \item $T$ is a quasi self-adjoint extensions of $\{H_{1,0},H_{2,0}\}$ if and only  if
  there exists $Q\subset \C^{2n}$ with $\dim Q=n$ such that
\begin{equation*}
     T= \{(y,[g])\in H_{1}:\; y(a)\in Q\}.
  \end{equation*}
\end{enumerate}}

\smallskip

The rest of this paper is organized as follows.
In Section 2, the basic concepts and useful results on the linear relations
are recalled in the first subsection.
In the second subsection, some important properties of the solutions of (1.1) and (\ref{HS2}) are presented.
The relationships between the minimal linear relations generated by (1.1) and (\ref{HS2})   are proved In subsection 3.1,
and characterizations of the minimal linear relations  are obtained in subsection 3.2.
These results extend  the corresponding results
from formally self-adjoint discrete linear Hamiltonian systems to formally non-self-adjoint ones.
In Section 4, we establish a bijective projection between all the quasi self-adjoint extensions
of $\{H_{1,0},H_{2,0}\}$ and all the self-adjoint extensions of the minimal linear relations $\mathbf{H_0}$ generated by $(\ref{HS3})$.
In Section 5, the case that $\I=[a,\infty)$ and $(\ref{HS4})$ is in the limit point case is concerned.
A complete characterization of  all the quasi self-adjoint extensions of $\{H_{1,0},H_{2,0}\}$ is obtained
by a  subspace $Q\subset \C^{2n}$ with $\dim Q=n$ in terms of boundary conditions.

 \section{  Preliminaries }

In this section, we first recall some the basic concepts and useful results on the linear relations.
In the second subsection,  some useful properties of solutions of  (\ref{HS1}) and (\ref{HS2})
are obtained.

\medskip

\subsection{Basic theory of  linear relations}

In this subsection, we  recall some basic results about linear relations.
The reader are referred to \cite{Arens1961,Cross1998,Shi2012,XuRen2024}.

Let $X$  be a complex Hilbert  space  with inner product $\langle\cdot, \cdot\rangle$.
Let $T$ and $S$ be  two linear relations (briefly, relations or named linear subspaces) in $X^2:=X\times X$.
By $\mathcal{LR}(X)$ denote the set of all linear relations of $X^2$. Denote
\begin{align*}
&\mathcal{D}(T):=\{y\in X:\; \exists\;(y,g)\in T\},\\
&\mathcal{R}(T):= \{g\in X: \;\exists\;(y,g)\in T \},\\
&\mathcal{N}(T):= \{y\in X: \;\exists\;(y,0)\in T \},\\
&T(y):=\{ g\in X:\; \exists\;(y, g)\in T\}.
\end{align*}
It can be easily verified that $T(0)=\{0\}$ if and only if $T$ can determine a unique linear operator from
$\mathcal{D}(T)$ into $X$ whose graph is just $T$.

Let $T, S\in \mathcal{LR}(X)$. The inverse of $T$, the adjoint of $T$, and the sum of $T$ and $S$ are defined as
\begin{align*}
&T^{-1}:= \{(g,y)\in X^2: \;(y,g)\in T \},\\
&T^*:=\{(y,g)\in X^2: \;\langle y, f\rangle=\langle g,x\rangle,\; {\forall}\;(x,f)\in T\},\\
&T+S:=\{(y,f+g)\in X^2: \; (y,f)\in T, \;  (y,g)\in S\}.
\end{align*}

\begin{definition}
Let $T\in \mathcal{LR}(X)$. Then
\begin{itemize}
\item [$(1)$] $T$ is said to be {\it Hermitian}  if $T\subset T^*$.
$T$ is said to be {\it symmetric}   if $T\subset T^*$ and $\mathcal{D}(T)$ is dense in $X$.
$T$ is said to be {\it self-adjoint}  if $T=T^*$.
\item [$(2)$]  $T$ is said to be {\it closed} if $T=\overline{T}$,
where $\overline{T}$ is the closure of $T$ in $X^2$.
By $\mathcal{CR}(X)$ denote the set of all the closed linear relations of $X^2$.
\item [$(3)$] Let $\ld\in \C$. The subspace $(\mathcal{R}(T-\ld))^\bot$ and the number
$d_{\ld}(T) := \dim (\mathcal{R}(T-\ld))^\bot$ are called the {\it deficiency
space} and {\it deficiency index} of $T$ with $\ld$, respectively.
\end{itemize}
\end{definition}

The following are some basic results on linear relations.

\begin{lemma} \cite{Arens1961,Shi2012}\label{Arens}
Let $T\in \mathcal{LR}(X)$. Then
\begin{itemize}
\item [$(1)$]$T^*$ is a closed subspace in $X^2$.
\item [$(2)$]$T^*=(\overline{T})^*$ and $T^{**}=\overline{T}$.
\item [$(3)$] $\mathcal{N}(T^*)=(\mathcal{R}(T))^{\bot}=(\mathcal{R}(\overline{T}))^{\bot}$.
\item [$(4)$] $d_{\ld}(T)=d_{\ld}(\overline{T})$ for all $\ld\in \C$.
If $T$ is Hermitian, then $d_{\ld}(T)$ is constant in the upper and lower half-planes,  respectively. Denote
\begin{equation*}
  d_{\pm}(T) := d_{\pm i}(T)
\end{equation*}
for an  Hermitian linear relation $T$, and call $d_{\pm}(T)$ the positive and negative deficiency indices of  $T$, respectively.
\end{itemize}
\end{lemma}

\begin{definition}\cite[Definitions 2.2,3.1]{XuRen2024}
\begin{enumerate}
\item  Let $T,S\in \mathcal{CR}(X)$. $\{T,S\}$ is  said to be a dual pair  if
\begin{equation*}
 \langle f,y\rangle=\langle x,g\rangle,\quad {\forall}\;(x,f)\in T,\;{\forall}\;(y,g)\in S,
\end{equation*}
or equivalently if $T\subset S^*$ ($\Leftrightarrow S\subset T^*$).
\item Let $K\in \mathcal{CR}(X)$. $K$ is called a proper extension of a dual pair $\{T,S\}$
if $T\subset K\subset S^*$.
The set of all proper extensions of  dual pair $\{T,S\}$ is denoted by {\rm Ext$\{T,S\}$}.
\item An extension $K\in {\rm Ext}\{T,S\}$ is called a quasi self-adjoint extension of $\{T,S\}$ if
\begin{equation*}
\dim(\mathcal{D}(K)/\mathcal{D}(T))=\dim(\mathcal{D}(K^*)/\mathcal{D}(S)).
\end{equation*}
\item  A dual pair $\{T,S\}$  is called a correct dual pair if it admits a quasi self-adjoint extension.
\end{enumerate}
\end{definition}

In 1961, Arens [1] introduced the following important
decomposition for a closed subspace $T$ in $X^2$:
\begin{equation*}
T=T_s\oplus T_\infty,
\end{equation*}
where
\begin{equation*}
T_\infty :=\{(0,g)\in X^2: (0,g)\in T\},\;\;
T_s :=T\ominus T_{\infty}.
\end{equation*}
$T_s$ and $T_\infty$ are called the operator and pure multi-valued parts of $T$, respectively.
It is evident that
\begin{equation*}
\mathcal{D}(T_s)=\mathcal{D}(T),\quad \mathcal{R}(T_{\infty})=T(0),\quad T_\infty =\{0\}\times T(0).
\end{equation*}

The following results are obtained.

\begin{lemma}\label{dual pair}\cite[Proposition 2.1, Theorems 3.1, 3.2]{XuRen2024}
Let  $T,S\in \mathcal{CR}(X)$ and $\{T,S\}$ be a dual pair with $(T^*)_s|_{\mathcal{D}(S)}=S_{s}$ and $(S^*)_s|_{\mathcal{D}(T)}=T_{s}$.
Suppose that $S^*(0)\cap \mathcal{N}(T^*)=\{0\}$ and  $T^*(0)\cap \mathcal{N}(S^*)=\{0\}$. Then
\begin{enumerate}
  \item
  \begin{equation*}
  \dim(\mathcal{D}(T^*)/\mathcal{D}(S))=\dim(\mathcal{D}(S^*)/\mathcal{D}(T)).
  \end{equation*}
  \item Let $K\in {\rm Ext}\{T,S\}$. Then
   \begin{align*}
  \dim(\mathcal{D}(T^*)/\mathcal{D}(S))&=\dim(\mathcal{D}(S^*)/\mathcal{D}(K))+\dim(\mathcal{D}(T^*)/\mathcal{D}(K^*))\nonumber\\
  &=\dim(\mathcal{D}(K)/\mathcal{D}(T))+\dim(\mathcal{D}(K^*)/\mathcal{D}(S)).
  \end{align*}
  \item Let $K$ be a quasi self-adjoint extension of $\{T,S\}$.
  Then $\dim(\mathcal{D}(T^*)/\mathcal{D}(S))$ is even and
   \begin{align*}
  \frac{1}{2}\dim(\mathcal{D}(T^*)/\mathcal{D}(S))&=\dim(\mathcal{D}(S^*)/\mathcal{D}(K))=\dim(\mathcal{D}(T^*)/\mathcal{D}(K^*))\nonumber\\
  &=\dim(\mathcal{D}(K)/\mathcal{D}(T))=\dim(\mathcal{D}(K^*)/\mathcal{D}(S)).
  \end{align*}
  \item Let $K\in {\rm Ext}\{T,S\}$ and suppose
  \begin{align*}
  \dim(\mathcal{D}(S^*)/\mathcal{D}(K))=\dim(\mathcal{D}(K)/\mathcal{D}(T)).
  \end{align*}
  Then $K$ is a quasi self-adjoint extension of $\{T,S\}$.
\end{enumerate}
\end{lemma}

Next, we introduce a form on $X^2\times X^2$ by
\begin{align*}
[(x,f):(y,g)]:=\langle f,y\rangle-\langle x,g\rangle,\quad
(x,f),(y,g)\in X^2.
\end{align*}
It can be easily verified that $[:]$ is a conjugate bilinear and
skew-Hermitian map from $X^2\times X^2$ into $\C$.

\begin{theorem}\label{closure T}
Let  $T,S\in \mathcal{LR}(X)$ and $\{T,S\}$ be a dual pair.
Then  $[T:S]=[T:T^*]=0$, and
\begin{align}\label{closure T1}
    \overline{T}=\{(x,f)\in S^*:\;[(x,f):T^*]=0\}.
\end{align}
\end{theorem}

\begin{proof}
Let $\{T,S\}$ be a dual pair.
Then it follows that $T\subset S^*$ and $S\subset T^*$, and  consequently $\overline{T}\subset S^*$.
It is clear that  $[T:S]=[T:T^*]=0$.

By $T_1$ denote the set on the right-hand side of (\ref{closure T1}).
Then $T_1$ is closed since $T^*$ is closed.
It is clear that $T\subset T_1$, and consequently  $\overline{T}\subset T_1$.

On the other hand, for any $(x,f)\in T_1$, it follows from the right-hand side of (\ref{closure T1})
that $[(x,f):(y,g)]=\langle f,y\rangle-\langle x,g\rangle=0$ for all $(y,g)\in T^*$.
This yields that $T_1\subset T^{**}=\overline{T}$, where (2) of  Lemma \ref{Arens} is used.
Hence, $\overline{T}=T_1$. The proof is complete.
\end{proof}

\medskip

\subsection{ Properties of solutions}

In this subsection, we discuss the properties of solutions of  (\ref{HS1}) and (\ref{HS2}), and as well as
the corresponding inhomogeneous systems.

First, we give some notations. Since $b$ may be finite or infinite,  we introduce the following
convention for briefness in the sequent discussion: $b+1$ means
$+\infty$ in the case of $b=+\infty$.
Denote
\begin{eqnarray*}
  l(\I):=\{y:\;y=\{y(t)\}_{t=a}^{b+1}\subset \C^{2n}\},
\end{eqnarray*}
and
\begin{equation*}
l^2_{W}(\I):=\left\{y\in l(\I):\;\sum_{t\in
\I}R(y)^*(t)W(t)R(y)(t)<+\infty\right\}
\end{equation*}
with the semi-scalar product
\begin{eqnarray*}
\langle y,z\rangle:=\sum_{t\in \I}R^*(z)(t)W(t)R(y)(t).
\end{eqnarray*}
Further,  define $\|y\|:=(\langle y,y\rangle)^{1/2}$ for $y\in l^2_{W}(\I)$.
Since the weight function $W(t)$ may be singular on $\I$,  we introduce the following quotient space:
\begin{equation*}
L^2_{W}(\I):=l^2_W(\I)/\{y\in l^2_W(\I):\; \|y\|=0\}.
\end{equation*}
By  \cite[Lemma 2.5]{Shi2006}, $L^2_W(\I)$ is a Hilbert space with
the inner product $\langle \cdot,\cdot\rangle$.
For a function $y\in l^2_W(\I)$, we denote by $[y]$   the corresponding
equivalent class in $L^2_W(\I)$.
Denote
\begin{eqnarray*}
&&l^2_{W,0}(\I):=\{y\in l^2_{W}(\I):\;
\exists\; s,k\in \I \;{\rm s.\; t.}\; y(t)=0
\;{\rm for}\; t\leq s\; {\rm and}\; t\ge k+1\}.
\end{eqnarray*}

The result on $\Ls_i$ defined by (1.1) and (\ref{HS2}) is obtained.

\begin{lemma}\label{liuville lemma}
For any $x, y\in l(\I)$ and any $[s,k]\subseteq\I$,
\begin{align}\label{liuville}
\sum_{t=s}^{k}[R(y)^*(t)(\Ls_2 x)(t)- (\Ls_1y)^*(t)R(x)(t)]=
y^*(t)Jx(t)|_{s}^{k+1}.
\end{align}
Further, let $y_i(t,\ld)$ be any solution  of  $(\Ls_iy)(t)=\ld W(t)R(y)(t)$ on $\I$, $i=1,2$.
Then there exists a constant $c$ such that
\begin{align}\label{liuville 2}
y^*_2(t,\bar{\ld})Jy_1(t,\ld)=c, \quad \forall\;t\in \I.
\end{align}
\end{lemma}

\begin{proof}
For any $x, y\in l(\I)$, it can be easily verified  that
\begin{align}\label{liuville Delta}
&R(y)^*(t)(\Ls_1 x)(t)- (\Ls_2y)^*(t)R(x)(t)\nonumber  \\\nonumber
=&R(y)^*(t)[J\De x(t)-P(t)R(x)(t)]-[(J\De y(t)-P^*(t)R(y)(t)]^*JR(x)(t)\\\nonumber
=&\De(y^*(t)Jx(t))+\De y^*(t)\left[\left( \begin{array}{cc}
                   I_n & 0 \\
                   0 & 0 \\
                 \end{array}
               \right)J+J\left(
                 \begin{array}{cc}
                   0 & 0 \\
                   0 & -I_n \\
                 \end{array}
               \right)\right]\De x(t)\\
 =&\De(y^*(t)Jx(t)),
\end{align}
where the fact
\begin{align*}
  R(y)(t)=y(t)+\left(
                 \begin{array}{cc}
                   I_n & 0 \\
                   0 & 0 \\
                 \end{array}
               \right)\De y(t)
  =y(t+1)+\left(
                 \begin{array}{cc}
                   0 & 0 \\
                   0 & -I_n \\
                 \end{array}
               \right)\De y(t)
\end{align*}
is used. Therefore, (\ref{liuville}) is followed directly by summing (\ref{liuville Delta}) from $s$ to $k$ on both side.

Let $y_i(t,\ld)$ be any solution  of  $(\Ls_iy)(t)=\ld W(t)R(y)(t)$ on $\I$. Then it follows that
\begin{align*}
[R(y_1)^*(t,\ld)(\Ls_2 y_2)(t,\bar{\ld})- (\Ls_1y_1)^*(t,\ld)R(y_2)(t,\bar{\ld})]=0,\quad t\in \I.
\end{align*}
This, together with (\ref{liuville Delta}),  yields $\De(y^*_1(t,\ld)Jy_2(t,\bar{\ld}))=0$ for all $t\in \I$.
Therefore, (\ref{liuville 2}) holds.
The proof is complete.
\end{proof}

For  any fixed $c_0\in \I$ and any $\ld\in \C$, by  $Y_i(t,\ld)$ denoted  the fundamental
solution matrix of  $(\Ls_iy)(t)=\ld W(t)R(y)(t)$ on $\I$ with $Y_i(c_0,\ld)=I_{2n}$ in the sequel.
It follows from (\ref{liuville 2}) that
\begin{align}\label{green 1 }
Y_2^*(t,\bar{\ld})JY_1(t,\ld)=Y_2^*(c_0,\bar{\ld})JY_1(c_0,\ld)=J, \quad \forall\;t\in \I.
\end{align}
In addition, for any  $\ld\in \C$, let $y(t,\ld)$   be any solution  of $(1.1)$ on $\I$.
Then it can be derived directly from   (\ref{HS blocked}) that
\begin{equation}\label{R(y)}
  R(y)(t,\ld)=\left(\begin{array}{cc}
                  (I_n-A(t))^{-1} & (I_n-A(t))^{-1}(B(t)+\ld W_2(t)) \\
                  0 & I_n
                \end{array}
                   \right)y(t,\ld),\quad \quad t\in \I.
\end{equation}

Next, we use the method of  variation of constants to obtain the solutions of  the following inhomogeneous systems
for $i=1,2$:
\begin{align}\label{IHS}
(\Ls_iy)(t)-\ld W(t)R(y)(t)=W(t)R(g)(t),\quad t\in \I.
\end{align}

\begin{lemma}\label{V C formulae}
Let $\ld\in \C$ and $i=1,2$. The solution $y_i(t,\ld)$ of $(\ref{IHS})$ with initial value $y_i(c_0,\ld)=y_0$ can be expressed uniquely by
\begin{align}\label{Const V}
y_i(t,\ld)=\left\{\ba{ll}
Y_i(t,\ld)y_0+Y_i(t,\ld)J\sum^{c_0-1}_{s=t}R(Y_{3-i})^*(s,\bar{\ld})W(s)R(g)(s),&a\le
t\le c_0-1,
\\
y_0,&t=c_0,\\
Y_i(t,\ld)y_0-Y_i(t,\ld)J\sum^{t-1}_{s=c_0}R(Y_{3-i})^*(s,\bar{\ld})W(s)R(g)(s),
&c_0+1\le t\le b+1.\ea\right.
\end{align}
\end{lemma}

\begin{proof} We show (\ref{Const V}) holds for $i=2$, and the other case $i=1$ can be proved similarly.

Let  $y(t)$ be a solution of ($\ref{IHS}$) with $i=2$. Here we omit the parameter $\ld$ for shortness.
By the method of  variation of constants, $y(t)$ can be written as $y(t)=Y_2(t,\ld)c(t)$.
It can be easily verified that
\begin{align*}
  &\Delta y(t)=\Delta Y_2(t,\ld)c(t+1)+Y_2(t,\ld)\Delta c(t),\\
  & R(y)(t)=R(Y_2)(t,\ld)c(t+1)-{\rm diag}\{0,I_n\}Y_2(t,\ld)\Delta c(t).
\end{align*}
So,  we can obtain  that
\begin{align}\label{J delta y 2}
&J\Delta y(t)-(P^*(t)+\ld W(t))R(y)(t)=J[\Delta Y_2(t,\ld)c(t+1)+Y_2(t,\ld)\Delta c(t)]\nonumber\\
&-(P^*(t)+\ld W(t))[R(Y_2)(t,\ld)c(t+1)-{\rm diag}\{0,I_n\}Y_2(t,\ld)\Delta c(t)]\nonumber\\
=&[J\Delta Y_2(t,\ld)-(P^*(t)+\ld W(t))R(Y_2)(t,\ld)]c(t+1)\\
&+[J+(P^*(t)+\ld W(t)){\rm diag}\{0,I_n\}]Y_2(t,\ld)\Delta c(t).\nonumber
\end{align}
By the definition of  $Y_2(t,\ld)$, it is clear that  $J\Delta Y_2(t,\ld)-(P^*(t)+\ld W(t))R(Y_2)(t,\ld)=0$.
In addition,  it can be easily verified that
\begin{align*}
J+(P^*(t)+\ld W(t)){\rm diag}\{0,I_n\}
=\left(
   \begin{array}{cc}
     0 & -(I_n-A^*(t)) \\
     I_n & B^*(t)+\ld W_2(t) \\
   \end{array}
 \right).
\end{align*}
Inserting these into (\ref{J delta y 2}), one can get that
\begin{align}\label{J delta y}
J\Delta y(t)-(P^*(t)+\ld W(t))R(y)(t)=\left(
   \begin{array}{cc}
     0 & -(I_n-A^*(t)) \\
     I_n & B^*(t)+\ld W_2(t) \\
   \end{array}
 \right)Y_2(t,\ld)\Delta c(t).
\end{align}
Thus, by inserting (\ref{J delta y}) into the left side of ($\ref{IHS}$) with $i=2$, we can obtain that
\begin{align}\label{delta c}
\left(
   \begin{array}{cc}
     0 & -(I_n-A^*(t)) \\
     I_n & B^*(t)+\ld W_2(t) \\
   \end{array}
 \right)Y_2(t,\ld)\Delta c(t)=W(t)R(g)(t).
\end{align}
By $(A_1), $ $I_n-A^*(t)$ is invertible on $\I$. So, the first matrix on the left side of
(\ref{delta c}) is invertible.
In addition, it follows from (\ref{green 1 }) that
\begin{equation*}
 Y_2^{-1}(t,\ld)=-JY_1^*(t,\bar{\ld})J.
\end{equation*}
So, we can get by (\ref{delta c}) that
\begin{align}\label{Delta c 1}
\Delta c(t)&=-JY_1^*(t,\bar{\ld})J\left(
   \begin{array}{cc}
     0 & -(I_n-A^*(t)) \\
     I_n & B(t)+\ld W_2(t) \\
   \end{array}
 \right)^{-1}W(t)R(g)(t)\nonumber\\
  &=-J\left[\left(
   \begin{array}{cc}
        (I_n-A(t))^{-1} & (I_n-A(t))^{-1}( B(t)+\bar{\ld} W_2(t))\\
         0& I_n
   \end{array}
 \right)Y_1(t,\bar{\ld})\right]^*W(t)R(g)(t).
\end{align}
It follows from (\ref{R(y)}) that
\begin{equation*}\label{R(Y)}
  R(Y_1)(t,\bar{\ld})=\left(\begin{array}{cc}
                  (I_n-A(t))^{-1} & (I_n-A(t))^{-1}(B(t)+\bar{\ld} W_2(t)) \\
                  0 & I_n
                \end{array}
                   \right)Y_1(t,\bar{\ld}),\quad \quad t\in \I.
\end{equation*}
Therefor,  (\ref{Delta c 1})  can be written as
\begin{align}\label{delta c 2}
\Delta c(t)=-JR(Y_1)^*(t,\bar{\ld})W(t)R(g)(t).
\end{align}
By summing (\ref{delta c 2}) from $c_0$ to $t-1$ on both side and by using  $c(c_0)=y_0$, one get that
\begin{align*}
 c(t)=y_0-J\sum_{s=c_0}^{t-1}R(Y_1)^*(s,\bar{\ld})W(s)R(g)(s)
\end{align*}
for $t> c_0$, and
\begin{align*}
 c(t)=y_0+J\sum_{s=t}^{c_0-1}R(Y_1)^*(s,\bar{\ld})W(s)R(g)(s)
\end{align*}
for $t< c_0$.
Inserting these into $y(t)=Y_2(t,\ld)c(t)$, one get (\ref{Const V}) holds for $i=2$.
The proof is complete.
\end{proof}

\begin{remark}
Lemmas \ref{liuville lemma} and \ref{V C formulae}  extend  the corresponding results
from formally self-adjoint discrete linear Hamiltonians \cite[Lemma 3.2]{Ren2014Laa} and \cite[Lemma 2.2]{Shi2006}
to formally non-self-adjoint ones.
\end{remark}

\section{The maximal and minimal linear relations}

In the first subsection, the relationship between the minimal linear relations generated by (1.1) and (\ref{HS2})   are established.
In the second subsection, the characterizations of the minimal linear relations  are obtained. 

\smallskip

\subsection{The  relationship of the maximal and minimal linear relations }

In this subsection, we discuss the relationship  between the minimal linear relations generated by (1.1) and (\ref{HS2}).

Denote
\begin{eqnarray*}
&&H_i:=\{([y],[g])\in (L^2_W(\I))^2: \;{\rm \exists\;}
y\in[y]\; {\rm s.t. }\; (\Ls_iy)(t)=W(t)R(g)(t),\; t\in\I\}, \\
&&H_{i,00}:=\{([y],[g])\in H: \;{\rm \exists\;}
y\in [y]\;{\rm s.t.}\;  y\in l^2_{W,0}(\I)\;{\rm and}\; (\Ls_iy)(t)=W(t)R(g)(t),\; t\in
\I\},\nonumber\\
&&H_{i,0}:=\overline{H}_{i,00},
\end{eqnarray*}
where $H_i$, $H_{i,00}$, and $H_{i,0}$  are called the maximal, pre-minimal,
and minimal relations corresponding to $\Ls_i$, separately.

Fix $c_0\in \I$ and let $\ld\in \C$. For $i=1,2$, still by $Y_i(t,\ld)$ denote the fundamental solution matrix of
 $(\Ls_iy)(t)=\ld W(t)R(y)(t)$ with $Y_i(c_0,\ld)=I_{2n}$.
For any finite subinterval $\I_1\subseteq \I$ with $c_0\in \I_1$, denote
\begin{align*}
\Phi_{i}(\ld, \I_1):=\sum_{t\in \I_1}R(Y_i)^*(t,\ld)W(t)R(Y_i)(t,\ld).
\end{align*}
It is evident that  $\Phi_{i}(\ld, \I_1)$, $i=1,2$, are  $2n\times 2n$ positive
semi-definite matrix,  and they dependent on $\ld$ and $\I_1$.

In the case that $\I$ is infinite, that is,  $\I=[a,\infty)$ or
$\I=(-\infty, b]$ or $\I=(-\infty,\infty)$, where $a,b$ are finite
integers, we introduce the following subspaces of  $L^2_{W}(\I)$,
respectively:
\begin{align*}
&L^2_{W,1}([a,\infty)):=\{[y]\in L^2_{W}(\I): \exists \; s\in \I\;
{\rm s.t.}\; W(t)R(y)(t)=0 \;{\rm for}\; t\geq s\},\\
&L^2_{W,1}((-\infty, b]):=\{[y]\in L^2_{W}(\I):\exists \; s\in \I\; {\rm
s.t.}\; W(t)R(y)(t)=0 \;{\rm for}\; t\leq s\},\\
&L^2_{W,1}((-\infty,\infty)):=\{[y]\in L^2_{W}(\I): \exists \; s \in
\I\;{\rm s.t.}\; W(t)R(y)(t)=0 \; {\rm for}\; |t|\geq |s|\}.
\end{align*}
It is evident  that $L^2_{W,1}(\I)$ is dense in $L^2_{W}(\I)$.

For $i=1,2$, map $\phi_{i,\ld,\I}$ is defined by
\begin{align*}
\phi_{i,\ld,\I}: L^2_W(\I)\to \C^{2n},\quad [g]\mapsto \sum_{t\in\I}R(Y_i)^*(t,\ld)W(t)R(g)(t)
\end{align*}
in the case that $\I$ is finite, and by
\begin{align*}
\phi_{i,\ld,\I}: L^2_{W,1}(\I)\to \C^{2n},\quad [g]\mapsto
\sum_{t\in\I}R(Y_i)^*(t,\ld)W(t)R(g)(t)
\end{align*}
in the case that $\I$ is infinite.

The following results are obtained.

\begin{lemma}\label{Phi I ld}  Let $i=1,2$.
\begin{itemize}
\item[\rm {(1)}] There exists a finite subinterval $\I_0:=[s_0,t_0]$ with $c_0\in \I_0\subseteq \I$  such that
\begin{align*}
{\rm rank}\,\Phi_{i}(\ld,\I_0)={\rm rank}\,\Phi_{i}(\ld,\I_1),\quad {\rm
Ran}\,\Phi_{i}(\ld,\I_0)={\rm Ran}\,\Phi_{i}(\ld,\I_1)
\end{align*}
for any finite subinterval $\I_1$ satisfying $\I_0\subseteq \I_1\subseteq \I$. In addition,
${\rm rank}\,\Phi_{i}(\ld,\I_0)$ and ${\rm Ran}\,\Phi_{i}(\ld,\I_0)$,  $i=1,2$, are both independent of $\ld$.
In the sequel, we  denote
\begin{align*}
{\rm rank}\,\Phi_{i}:={\rm rank}\,\Phi_{i}(\ld,\I_0),\quad {\rm Ran}\,\Phi_{i}:={\rm Ran}\,\Phi_{i}(\ld,\I_0).
\end{align*}
\item[\rm {(2)}] For all $\ld\in\C$, ${\rm Ran}\,\phi_{i,\ld,\I}={\rm Ran}\,\Phi_{i}$,
and consequently  ${\rm Ran}\,\phi_{i,\ld,\I}$ is independent of $\ld$.
\item[\rm {(3)}] In the case that $\I$ is finite,
\begin{align*}
L^2_W(\I)={\rm Ker}\,\phi_{i,\ld,\I}\oplus\{[Y_i(\cdot,\ld)\xi]:\;\xi\in {\rm Ran}\,\Phi_{i}\}.
\end{align*}
\item[\rm {(4)}] In the case that $\I$ is infinite, let $r_i={\rm rank}\, \Phi_{i}$.
Then there exist  linearly independent elements $[\af_{i,j}]\in L^2_{W,1}(\I)$, $1\le j\le r_i$, such that
\begin{align}\label{L^2_{W,1}}
L^2_{W,1}(\I)={\rm Ker}\,\phi_{i,\ld,\I}\dotplus {\rm
span}\{[\af_{i,1}],[\af_{i,2}],\ldots,[\af_{i,r_i}]\}.
\end{align}
\item[\rm {(5)}]  For all  $\ld\in\C$,
\begin{equation}\label{Ker phi}
 {\rm Ker}\, \phi_{i,\ld,\I} \subset {\rm Ran}\,(H_{3-i,00}-\bar{\ld} ).
\end{equation}
\end{itemize}
\end{lemma}

\begin{proof}
Since the proofs of  \cite[Lemmas 3.1, 3.2]{Ren2011} and the proofs of (1) and (2) of \cite[Lemmas 3.3]{Ren2011}
do not require $P^*=P$,
the  assertions (1)-(4) can be proved  with  similar arguments. We omit the details.

Next, we show assertion (5) holds.
We only need to show (\ref{Ker phi}) holds for $i=1$,
and the other can be proved similarly.
The proof is divided into two parts.

{\bf Case 1.}  $\I=[a,b]$ is finite.
For any $\ld\in \C$ and any  $[g]\in{\rm Ker}\,\phi_{1,\ld,\I}$, let $y$ be a solution of the initial valued problem
\begin{equation*}
  (\Ls_2y)(t)=\bar{\ld} W(t)R(y)(t)+W(t)R(g)(t),\quad y(a)=0,\quad t\in\I.
\end{equation*}
Then it follows from (\ref{Const V}) that
\begin{align}\label{y(t)}
y(t)=-Y_2(t,\bar{\ld})J\sum^{t-1}_{s=a}R(Y_1)^*(s,\ld)W(s)R(g)(s), \quad a+1 \le
t\le b+1.
\end{align}
Since  $[g]\in{\rm Ker}\, \phi_{1,\ld,\I} $,  it follows  from (\ref{y(t)}) that
\begin{align*}
y(b+1)&=-Y_2(b+1,\bar{\ld})J\sum_{s\in\I}R(Y_1)^*(s,\ld)W(s)R(g)(s)=-Y_2(b+1,\bar{\ld})J\phi_{1,\ld,\I}([g])=0.
\end{align*}
This yields that $([y],[g])\in H_{2,00}-\bar{\ld}$ and $[g]\in {\rm Ran}\,(H_{2,00}-\bar{\ld})$.
By the arbitrariness of $[g]$, one has that ${\rm Ker}\, \phi_{1,\ld,\I}\subset {\rm Ran}\, (H_{2,00}-\bar{\ld})$.

{\bf Case 2.}  $\I$ is infinite. We only show that the assertions hold in the case that
$\I=[a,\infty)$. For the other two cases, it can be proved with similar arguments.

Fix $\ld\in\C$.
For any  $[g]\in{\rm Ker}\,\phi_{1,\ld,\I}\subseteq L^2_{W,1}(\I)$, there exists an integer $k\ge t_0+1$ such that
$W(t)R(g)(t)=0$ for all $t\geq k$.
Let $y$ be a solution of the initial valued problem
\begin{equation*}
  (\Ls_2y)(t)=\bar{\ld} W(t)R(y)(t)+W(t)R(g)(t),\quad y(k)=0,\quad t\in\I.
\end{equation*}
Then it follows from (\ref{Const V}) that
\begin{align}\label{y(t)2}
y(t)=Y_2(t,\bar{\ld})J\sum^{k-1}_{s=t}R(Y_1)^*(s,\ld)W(s)R(g)(s), \quad a\leq t\leq k-1,
\end{align}
and $y(t)=0$ for $t\ge k$. It is clear  that $y\in l_{W}^2(\I)$.
Since $[g]\in{\rm Ker}\, \phi_{1,\ld,\I}$ and $W(t)R(g)(t)=0$ for $t\geq k$,
we get by (\ref{y(t)2}) that
\begin{align*}
y(a)&=Y_2(a,\bar{\ld})J\sum^{k-1}_{s=a}R(Y_1)^*(s,\ld)W(s)R(g)(s)=Y_2(a,\bar{\ld})J\sum^{\infty}_{s=a}R(Y_1)^*(s,\ld)W(s)R(g)(s)\\
&=Y_2(a,\bar{\ld})J\phi_{1,\ld,\I}([g])=0.
\end{align*}
This yields that  $([y],[g])\in H_{2,00}-\bar{\ld}I $ and $[g]\in {\rm Ran}\,(H_{2,00}-\bar{\ld}I)$.
Hence,   ${\rm Ker}\, \phi_{1,\ld,\I}\subset {\rm Ran}\, (H_{2,00}-\bar{\ld}I)$.
The proof is complete.
\end{proof}

Based on Lemma \ref{Phi I ld}, we obtained
the relationship between the minimal linear relations generated by (1.1) and (\ref{HS2}).

\begin{theorem}\label{T0*=TM+}
Assume that $(A_1)$ holds. Then $H_2=H_{1,0}^*$,  $H_1=H_{2,0}^*$,
and consequently  $\{H_{1,0},H_{2,0}\}$ is a dual pair.
\end{theorem}

\begin{proof}
The main idea of the  proof is  similar as that of \cite[Theorem 3.1]{Ren2011}, which require $P^*=P$.
We only prove $H_1=H_{2,0}^*$ and the other can be showed similarly.
Since $H_{2,0}^*=H_{2,00}^*$, it suffices to show  $H_1=H_{2,00}^*$,

First, we show  $H_1\subseteq H_{2,00}^*$. Let $([y],[g])\in H_1$. Then there exists $y\in [y]$ such that
\begin{equation*}
  (\Ls_1y)(t)=W(t)R(g)(t), \quad t\in \I.
\end{equation*}
For any $ ([z],[h])\in H_{2,00}$, there exists $z\in [z]$ such that
\begin{equation*}
  (\Ls_2 z)(t)=W(t)R(h)(t), \quad t\in \I.
\end{equation*}
Then, by using  (\ref{liuville}) we get that
\begin{align*}
&\langle [h],[y]\rangle-\langle [z],[g]\rangle =
\sum_{t\in\I}\left[ R(y)^*(t)W(t)R(h)(t)-  R(g)^*(t)W(t)R(z)(t)\right]\\
=&\sum_{t\in\I}\left[ R(y)^*(t)(\Ls_2 z)(t)-(\Ls_1y)^*(t)R(z)(t) \right]=(y^*Jz)(t)|_a^{b+1}=0,
\end{align*}
where $z\in l^2_{W,0}(\I)$ is used. Hence,  $H_1\subseteq H_{2,00}^*$.

Next, we show $H_{2,00}^*\subseteq H_1$.  Let $([y], [g])\in H_{2,00}^*$.
Then  for any $([z], [h])\in H_{2,00}$, there exists $z\in[z]$ with $z\in l^2_{W,0}(\I)$ such that
\begin{equation*}
  (\Ls_2z)(t)=W(t)R(h)(t), \quad t\in \I.
\end{equation*}
and  it follows that
\begin{align}\label{liuville 3}
0=\langle [h],[y]\rangle-\langle [z],[g]\rangle =\sum_{t\in\I}[R(y)^*(t)W(t)R(h)(t)-R(g)^*(t)W(t)R(z)(t)].
\end{align}
For the above given $[g]$, let  $y_1\in l(\I)$ be a solution of the following equation
\begin{equation*}
  (\Ls_1 y_1)(t)=W(t)R(g)(t), \quad t\in \I.
\end{equation*}
The existence of $y_1$ is guaranteed by Assumption $(A_1)$.
Since  $y_1$ may not belong to $[y]$, we need to find some $y_0\in [y]$
such that \begin{equation*}
  (\Ls_1 y_0)(t)=W(t)R(g)(t), \quad t\in \I.
\end{equation*}
Again by using  (\ref{liuville}), one can get  that
\begin{align}\label{liuville 4}
&\sum_{t\in\I}[R(y_1)^*(t)W(t)R(h)(t)-R(g)^*(t)W(t)R(z)(t)]\nonumber\\
=&\sum_{t\in\I}[R(y_1)^*(t)(\Ls_2 z)(t)- (\Ls_1y_1)^*(t)R(z)(t)]=(y_1^*J z)(t)|_{a}^{b+1}=0,
\end{align}
where $z\in l^2_{W,0}(\I)$ is used.
By Combining (\ref{liuville 3}) and (\ref{liuville 4}), one has that for all $([z], [h])\in H_{2,00}$,
\begin{align*}
\sum_{t\in\I}R(y-y_1)^*(t)W(t)R(h)(t)]=0.
\end{align*}
Now we take $\ld=0$.  By (5) of   Lemma \ref{Phi I ld}, we get  that for any $[g]\in {\rm Ker}\, \phi_{1,0,\I}$
\begin{align}\label{liuville 5}
\sum_{t\in\I} R(y-y_1)^*(t)W(t)R(g)(t)=0.
\end{align}
On the other hand,  it follows from  $[g]\in {\rm Ker}\, \phi_{1,0,\I}$  that
\begin{align}\label{liuville 6}
\sum_{t\in\I}R(Y_1)^*(t,0)W(t)R(g)(t)=0.
\end{align}
So, by combining (\ref{liuville 5}) and (\ref{liuville 6}), we obtain  that
\begin{align}\label{sum 0}
\sum_{t\in\I}R(g)^*(t)W(t)R(y-y_1-Y_1(\cdot,0)\xi)(t)=0
\end{align}
holds for any $[g]\in {\rm Ker}\, \phi_{1,0,\I}$ and any $\xi\in \C^{2n}$.

In the case that $\I$ is finite, $[y_1], [y]\in L_W^2(\I)$ is evident.
Then  it follows from (3) of Lemma \ref{Phi I ld} that
there exists $\xi\in {\rm Ran}\,\Phi_{1}$ such that
\begin{align*}
[y-y_1-Y_1(\cdot,0)\xi]\in {\rm Ker}\, \phi_{1,0,\I}.
\end{align*}
This, together with (\ref{sum 0}), implies that $y-y_1-Y_1(\cdot,0)\xi=0$ in $ L_W^2(\I)$.
Set
\begin{equation*}
  y_0(t):=y_1(t)+Y_1(t,0)\xi,\quad t\in \I.
\end{equation*}
Then $y_0\in [y]$  and it satisfies
\begin{align*}
\Ls_1(y_0)(t)=(\Ls_1 y_1)(t)+(\Ls_1Y_1)(t,0)\xi=W(t)R(g)(t),\quad t\in \I.
\end{align*}

In the case that $\I$ is infinite, we only consider the case that $\I=[a,\infty)$.
For the other two cases, it can be proved with similar arguments.

It follows from (4) of Lemma \ref{Phi I ld} that there exist $[\af_{1,j}]\in L_{W,1}^2(\I)$, $1\le j\le r_1$,
such that (\ref{L^2_{W,1}}) holds.  Since each $[\af_{1,j}]\in L_{W,1}^2(\I)$,
there exists $k\ge t_0$, where $t_0$ is specified by (1) of Lemma \ref{Phi I ld}, such that
\begin{align*}
W(t)R(\af_{1,j})(t)=0, \quad t\ge k, \quad  1\le j\le r_1.
\end{align*}
Set $\I_1=[a,k]$. By $y|_{\I_1}$ denote the restriction of $y$ on $[a,k+1]$.
By (3) of Lemma \ref{Phi I ld}, there exist linearly independent elements
$\xi_1,\xi_2,\ldots,\xi_{r_1}\in {\rm Ran}\,\Phi_{1}$ such that
\begin{eqnarray*}
[(\af_j-Y_1(\cdot,0)\xi_j)|_{\I_1}]\in{\rm Ker}\, \phi_{1,0,\I_1}, \quad 1\le j\le r_1.
\end{eqnarray*}
Take $\bt_j(t)$, $1\le j\le r_1$, satisfying
\begin{align}\label{btj}
R(\bt_j)(t)=R(Y)(t,0)\xi_j\; {\rm for}\; t\in [a, k],\quad R(\bt_j)(t)=0
\;{\rm for}\; t\ge k+1.
\end{align}
It is clear that $\bt_1,\bt_2,\ldots,\bt_{r_1}$ are linearly independent elements in $L_{W,1}^2(\I)$ and
\begin{align}\label{btj 1}
[\bt_j|_{\I_1}]\in ({\rm Ker}\, \phi_{1,0,\I_1})^\bot,\quad
[\bt_j-\af_{1,j}]\in {\rm Ker}\, \phi_{1,0,\I},\quad 1\le j\le r_1.
\end{align}
Hence, it follows from (\ref{L^2_{W,1}}) that
\begin{align}\label{btj 2}
L^2_{W,1}(\I)={\rm Ker}\,\phi_{1,0,\I}\dotplus {\rm span}\{[\bt_1],[\bt_2],\ldots,[\bt_{r_1}]\}.
\end{align}

It is obvious that $[(y-y_1)|_{\I_1})\in L^2_{W}(\I_1)$. By (3) of Lemma \ref{Phi I ld}, one can get that
 there exists a $\xi\in {\rm Ran}\,\Phi_{\I_1}$ such that
\begin{align}\label{btj 3}
[(y-y_1-Y_1(\cdot,0)\xi)|_{\I_1}]\in {\rm Ker}\,\phi_{1,0,\I_1}.
\end{align}
By combining (\ref{btj}), (\ref{btj 1}) and (\ref{btj 3}),  we can obtain that
\begin{align}\label{btj 4}
&\sum_{t\in\I}R(\bt_j)^*(t)W(t)R(y-y_1-Y_1(\cdot,0)\xi)(t)\nonumber\\
=&\sum_{t\in\I_1}R(\bt_j)^*(t)W(t)R(y-y_1-Y_1(\cdot,0)\xi)(t)=0.
\end{align}
Hence, for any $[z] \in L^2_{W,1}(\I)$,  it follows from (\ref{sum 0}), (\ref{btj 2}) and (\ref{btj 4}) that
\begin{align*}
\sum_{t\in\I}R(z)^*(t)W(t) R(y-y_1-Y_1(\cdot,0)\xi)(t)=0.
\end{align*}
By the density of $L^2_{W,1}(\I)$ in $L^2_{W}(\I)$,  it follows that $[y-y_1-Y_1(\cdot,0)\xi]=0$.
Set $y_0(t):=y_1(t)+Y(t,0)\xi$. Then  $y_0\in [y]$ and satisfies
\begin{equation*}
 (\Ls_1y_0)(t)=(\Ls_1y_1)(t)=W(t)R(g)(t).
\end{equation*}
Therefore,  $([y],[g])\in H_1$.
By the arbitrariness  of $([y],[g])$ one has $H_{2,00}^*\subset H_1$.
The whole proof is complete.
\end{proof}

\subsection{ The characterizations of the  minimal linear relations}

In this subsection, we pay attention to the characterizations of the  minimal linear relations.

First we introduce the  definiteness  condition for  (1.1) and (\ref{HS2}):
\begin{itemize}
\item [\rm {($\bf A_{2}$)}]  There exists a finite subinterval
$\I_1\subseteq \I$, without loss of generality, we still denote $\I_1=\I_0=[s_0,t_0]$, such that
for any $g_i$  and $y_i$ with $y_i(k_i)\neq 0$ for some $k_i\in \I$, which satisfy $(\Ls_iy_i)(t)=W(t)R(g_i)(t)$ on $\I$,   $i=1,2$,
it follows that
\begin{align*}
\sum_{t\in \I_0}R(y_i)^*(t)W(t)R(y_i)(t)>0.
\end{align*}
\end{itemize}

\begin{remark}\label{remark}
\begin{enumerate}
  \item Assumption $(A_2)$ guarantees that for any
$\ld\in\C$,  every non-trivial solution $y_i(t)$ of $(\Ls_iy_i)(t)=\ld W(t)R(y_i)(t)$ on $\I$, satisfies
\begin{align*}
\sum_{t\in \I_0}R(y_i)^*(t)W(t)R(y_i)(t)>0.
\end{align*}
  \item Assume that $(A_1)$ and $(A_2)$ holds. With a similar argument as that of \cite[Theorem 4.2]{Ren2011}
one can show that   for any $([y_i],[g_i])\in H_i$,
there exists a unique $y_i\in [y_i]$ such that $(\mathscr{L}_iy_i)(t)=W(t)R(g_i)(t)$ on $\I$, $i=1,2$.
Therefore, we can just write  $(y_i,[g_i])\in H_i$ instead of $([y_i],[g_i])\in H_i$ under  $(A_1)$ and $(A_2)$ in the sequel.
\end{enumerate}

\end{remark}

Further, we can obtain the following result.

\begin{lemma}\label{H1H2} Assume that  $(A_1)$ and $(A_2)$ hold. Then
\begin{align*}
   \langle g_1,y_2\rangle-\langle y_1,g_2\rangle=(y_2^*Jy_1)(t)|_{a}^{b+1}.
\end{align*}
holds for any $ (y_i,[g_i])\in H_i$, $i=1,2$.
\end{lemma}

\begin{proof} Let $ (y_i,[g_i])\in H_i$,  $i=1,2$. Then it follows that
\begin{equation*}
  (\Ls_iy_i)(t)= W(t)R(g_i)(t), \quad t\in \I.
\end{equation*}
By using  (\ref{liuville}), one get that
\begin{align*}
   \langle g_1,y_2\rangle-\langle y_1,g_2\rangle &=\sum_{t\in \I}R(y_2)^*(t)W(t)R(g_1)(t)-R(g_2)^*(t)W(t)R(y_1)(t)\\
   &=\sum_{t\in\I} R(y_2)^*(t) (\Ls_1y_1)(t)-(\Ls_2y_2)^*(t)R(y_1)(t)=(y_2^*Jy_1)(t)|_{a}^{b+1}.
\end{align*}
Note that since each $ (y_i,[g_i])\in H_i$,  the limit $\lim_{t\to \infty}(y_2^*Jy_1)(t)$
exists finitely when $b=\infty$ and  $\lim_{t\to -\infty}(y_2^*Jy_1)(t)$ exists finitely when $a=-\infty$.
The proof is complete.
\end{proof}

In addition, the following result holds.

\begin{lemma}\label{patch lemma}
Assume that $(A_1)$ and $(A_{2})$
hold. Then for any given finite subset $\I_1=[s,k]$ satisfying
$\I_0\subseteq\I_1\subseteq \I$ and for any given $\af, \bt\in \C^{2n}$,
there exists $g_i=\{g_i(t)\}_{t=s}^{k+1}\subset \C^{2n}$ such that
the following boundary value problem:
\begin{align}\label{bvp}
\ba{l}(\Ls_i y)(t)=W(t)R(g_i)(t),\quad t\in \I_1,\\
y(s)=\af,\quad y(k+1)=\bt,\ea
\end{align}
has a solution $y_i=\{y_i(t)\}_{t=s}^{k+1}\subset \C^{2n}$, $i=1,2$.
\end{lemma}

\begin{proof}
The argument  is similar as that of \cite[Lemma 3.3]{Ren2014Laa} which require $P^*=P$,
We give the proof for completeness.
We only proof the result holds for $i=1$, and the other case can be proved with a similar argument.

Set
\begin{align*}
\langle x,y\rangle_1:=\sum_{t=s}^{k}R(y)^*(t)W(t)R(x)(t),
\end{align*}
and let $\phi_j$, $1\leq j\leq 2n$, be the linearly independent solutions of system $(\Ls_1\phi)(t)=0$ on $\I$.
Then we have
\begin{align}\label{ rank phi_i,j}
{\rm rank}\,(\langle \phi_i,\phi_j\rangle_1)_{1\leq i,j \leq 2n}=2n.
\end{align}
In fact, the linear algebraic system
\begin{align}\label{ inner phi_i,j}
(\langle \phi_i,\phi_j\rangle_1)_{1\leq i,j \leq 2n}C=0,
\end{align}
where $C=(c_1,c_2,\ldots,c_{2n})^T\in \C^{2n}$, can be written as
\begin{equation*}
  \langle \phi_i,\sum_{j=1}^{2n}\bar{c}_j\phi_j\rangle_1=0, \quad 1\leq i \leq 2n,
\end{equation*}
which yields that
\begin{equation*}
\langle
\sum_{j=1}^{2n}\bar{c}_j\phi_j,\sum_{j=1}^{2n}\bar{c}_j\phi_j\rangle_1 =0.
\end{equation*}
Since $\sum_{j=1}^{2n}\bar{c}_j\phi_j$ is still a solution of  $(\Ls_1\phi)(t)=0$ on $\I$,
it follows from $(A_2)$ that $\sum_{j=1}^{2n}\bar{c}_j\phi_j=0$.
This, together with the independence of $\phi_j$, $1\leq j\leq 2n$, yields that $C=0$.
This means that  (\ref{ inner phi_i,j}) has only a zero solution.
Therefore, (\ref{ rank phi_i,j}) holds.

Let $\af,\bt$ be any two given vectors in $\C^{2n}$.
By (\ref{ inner phi_i,j}), the linear algebraic system
\begin{align}\label{phi_i,j eq}
(\langle \phi_i,\phi_j\rangle_1)^T_{1\leq i,j\leq 2n}
C=(\phi_1(k+1),\phi_2(k+1), \ldots, \phi_{2n}(k+1))^*J\bt
\end{align}
has a unique solution $C_1\in \C^{2n}$.
Set $\psi_1=(\phi_1,\phi_2, \ldots, \phi_{2n})C_1$.
It follows from  (\ref{phi_i,j eq}) that
\begin{align}\label{psi_1,phi_i}
\langle \psi_1,\phi_i\rangle_1=\phi_i^*(k+1)J\bt,\quad 1\leq i\leq 2n.
\end{align}

Let $u(t)$ be a solution of the following  initial value problem:
\begin{align}\label{Ls_2x}
(\Ls_2x)(t)=W(t)R(\psi_1)(t),\quad t\in [s,k], \quad x(s)=0.
\end{align}
Since $(\Ls_1\phi_i)(t)=0$ for $t\in \I$ and $1\leq i\leq 2n$, we get by (\ref{Ls_2x}) and (\ref{liuville}) that
\begin{align}\label{psi_1,phi_i 2}
\langle
\psi_1,\phi_i\rangle_1&=\sum_{t=s}^k[R(\phi_i)^*(t)\Ls_2(u)(t)-\Ls_1(\phi_i)^*(t)R(u)(t)]\nonumber\\
&=(u,\phi_i)(t)|_{s}^{k+1}=\phi_i^*(k+1)Ju(k+1), \;\;1\leq i\leq 2n.
\end{align}
Since $\phi_1,\phi_2,\ldots,\phi_{2n}$ are linearly independent in $l(\I)$,
it follows from (\ref{psi_1,phi_i}) and (\ref{psi_1,phi_i 2}) that $u(k+1)=\bt$.
So, $u(t)$ is a solution of the following boundary value problem:
\begin{align*}
&(\Ls_2x)(t)=W(t)R(\psi_1)(t),\quad t\in [s,k],\\
&x(s)=0,\quad x(k+1)=\bt.
\end{align*}

On the other hand, the linear algebraic system
\begin{align}\label{phi_i,phi_j}
(\langle \phi_i,\phi_j\rangle_1)^T_{1\leq i,j\leq 2n}C=(\phi_1(s),\phi_2(s), \ldots, \phi_{2n}(s))^*J\af
\end{align}
has a unique solution $C_2\in\C^{2n}$ by (\ref{ rank phi_i,j}).
Set $\psi_2=(\phi_1,\phi_2, \ldots, \phi_{2n})C_2$. Then, by (\ref{phi_i,phi_j}) one can get that
\begin{align}\label{psi_2,phi_i}
\langle \psi_2,\phi_i\rangle_1=\phi_i^*(s)J\af,\quad 1\leq i\leq 2n.
\end{align}
Let $v(t)$ be a solution of the following  initial value problem:
\begin{align}\label{Ls_2x=-W}
\Ls_2(x)(t)=-W(t)R(\psi_2)(t),\quad x(k+1)=0,\quad t\in [s,k].
\end{align}
Since $(\Ls_2\phi_i)(t)=0$  for $t\in \I$ and $1\leq i\leq 2n$, we get by (\ref{liuville}) and (\ref{Ls_2x=-W}) that
\begin{equation*}
\langle \psi_2,\phi_i\rangle_1=\phi_i^*(s)Jv(s),\quad 1\leq i\leq 2n,
\end{equation*}
which, together with (\ref{psi_2,phi_i}), implies that $v(s)=\af$. So, $v(t)$ is
a solution of the following boundary value problem:
\begin{align*}
&(\Ls_2x)(t)=-W(t)R(\psi_2)(t),\quad t\in [s,k],\\
&x(s)=\af,\quad x(k+1)=0.
\end{align*}
Set $g_1=\psi_1-\psi_2$ and $y_1=u+v$.
Then $y_1$ is a solution of the boundary value problem (\ref{bvp}) with $i=1$.
This completes the proof.
\end{proof}

\begin{remark}\label{remark patch}
Lemma $\ref{patch lemma}$  is a extension of  the Patch lemma from the formally selfadjoint case \cite[Lemma 3.3]{Ren2014Laa}
to the formally non-selfadjoint case.
Based on this lemma, any two elements of $H_i$ can be patched up to construct another new element of $H_i$.
For example, take any $\I_1=[s, k]\supseteq\I_0$.
For any $ (x,[f]), (y,[g])\in H_1$, there exists $(z,[h])\in H_1$ satisfying
\begin{align*}
z(t)=x(t) \;\;{\rm for}\;\; t\leq s,\quad  z(t)=y(t)\;\;{\rm
for}\; \;t \ge k+1.
\end{align*}
In particular,  there exist $(z_{i,j},[h_{i,j}])\in H_{i}$ for $i=1,2$, $j=1,2,\ldots,2n$, satisfying
\begin{align*}
&z_{i,j}(s)=e_j;\quad z_{i,j}(t)=0,\quad t\ge k+1,
\end{align*}
where
\begin{align*}
    e_j:=(\underbrace{0,\ldots,0}_{j-1}, 1,0,\ldots,0)^T\in \C^{2n},
\quad 1\le j\le 2n.
\end{align*}
\end{remark}

\begin{corollary}\label{coro patch}
Assume that $(A_1)$ and $(A_2)$ hold. Then for  any given finite subset $\I_1=[s,k]$ satisying
$\I_0\subseteq\I_1\subseteq \I$  and any $\xi\in\C^{2n}$, there exists
$(y_i,[g_i])\in H_{i}$, $i=1,2$, satisfying
\begin{align*}
&y_i(s)=\xi;\quad y_i(t)=0,\quad t\ge k+1.
\end{align*}
\end{corollary}

Based on Lemma \ref{patch lemma} and Remark \ref{remark patch}, one can get the  following result.

\begin{theorem}\label{char H0}
If  $(A_1)$ and $(A_{2})$ hold, then for $i=1,2$,
\begin{align*}
H_{i,0}=\{(y,[g])\in H_i: (y^*Jx)(b+1)=(y^*Jx)(a)=0, \; \forall\; x\in \mathcal{D}(H_{3-i})\}.
\end{align*}
In particular,
\begin{align*}
H_{i,0}=\{(y,[g])\in H_i: y(a)=0, (y^*Jx)(b+1)=0, \; \forall\; x\in \mathcal{D}(H_{3-i})\}
\end{align*}
if $a$ is finite;
\begin{align*}
H_{i,0}=\{(y,[g])\in H_i: y(b+1)=0, (y^*Jx)(a)=0, \; \forall\; x\in \mathcal{D}(H_{3-i})\}
\end{align*}
if $b$ is finite;
\begin{align*}
H_{i,00}=H_{i,0}=\{(y,[g])\in H_i: y(a)=y(b+1)=0\}
\end{align*}
if both $a$ and $b$ are finite.
\end{theorem}

\begin{proof}\label{1}
For any $(x,[f])\in H_1, (y,[g])\in H_2$, we have by using (\ref{liuville}) that
\begin{align*}
\langle
f,y\rangle-\langle x,g\rangle=\sum_{t\in \I}[R(y)^*(t)(\Ls_1x)(t)- (\Ls_2y)^*(t)R(x)(t)]= (y^*Jx)(t)|_{a}^{b+1}.
\end{align*}
Therefore, by  Theorems \ref{closure T} and \ref{T0*=TM+}, one can get that
\begin{align}\label{H_{i,0}}
H_{i,0}=\{(y,[g])\in H_i: (y^*Jx)(b+1)=(y^*Jx)(a), \; \forall\; x\in \mathcal{D}(H_{3-i})\}.
\end{align}
For convenience,  denote
\begin{align*}
H^0_i=\{(y,[g])\in H_i: (y^*Jx)(b+1)=(y^*Jx)(a)=0, \; \forall\; y\in \mathcal{D}(H_{3-i})\}.
\end{align*}
Clearly, $H^0_i\subset H_{i,0}$.

We now show  $H_{i,0}\subset H^0_i$.
Fix any $(y,[g])\in H_{i,0}$. It follows from (\ref{H_{i,0}}) that for all $(x,[f])\in H_{3-i}$,
\begin{align}\label{(x,y)}
(y^*Jx)(b+1)= (y^*Jx)(a).
\end{align}
For any given $(x,[f])\in H_{3-i}$, by the discussion in Remark \ref{remark patch}, there exists $(z,[h])\in H_{3-i}$ such that
\begin{align*}
z(t)=0, \quad t\leq s_0;\quad z(t)=x(t), \quad t\ge t_0+1.
\end{align*}
Thus, it follows from (\ref{(x,y)}) that
\begin{equation*}
 (y^*Jx)(b+1)=(y^*Jz)(b+1)=(y^*Jz)(a)=0,
\end{equation*}
and consequently $(y^*Jx)(a)=0$ for all $(x,[f])\in H_{3-i}$.

In the case that $a$ is finite, there exist
$(z_j,[h_j])\in H$ with $z_j(a)=e_j$  for $j=1,2,\ldots,2n$ by Remark \ref{remark patch}.
It follows from $(y^*Jz_j)(a)=0$ for $j=1,2,\ldots,2n$  that $y(a)=0$.
With a similar argument, one can show that $y(b+1)=0$ if  $b$ is finite.
Therefore,  $H_{i,0}=H_{i,00}$ in the case that $\I=[a,b]$.
This completes the proof.
\end{proof}

\section{Quasi self-adjoint extensions  in general case}

In this section, we discuss the quasi self-adjointness of $H_{1,0}$ and $H_{2,0}$ in the general case.
To do this, we convert (1.1) and (\ref{HS2}) into a formally self-adjoint system by a homeomorphism.

Denote
\begin{align*}
  &E_1=\left(
     \begin{array}{cccc}
       0 & 0 & I_n & 0 \\
       I_n & 0 & 0 & 0 \\
       0 & I_n & 0 & 0 \\
       0 & 0 & 0 & I_n \\
     \end{array}
   \right),\quad
      E_2=\left(
     \begin{array}{cccc}
        I_n &0 & 0 & 0 \\
        0 & 0&I_n & 0 \\
       0 & 0 &0& I_n\\
       0 & I_n &0 &0
     \end{array}
   \right).
\end{align*}
It is clear that both $E_j$, $j=1,2$, are $4n\times 4n$ invertible matrices.

Denote
\begin{eqnarray*}
\mathbf{l}(\I):=\{\mathbf{y}:\;\mathbf{y}=\{\mathbf{y}(t)\}_{t=a}^{b+1}\subset \C^{4n}\},
\end{eqnarray*}
For any two $y_i\in l(\I)$,  $i=1,2$, we denote
\begin{equation}\label{y y y}
\mathbf{y}:=E_1\left(
                 \begin{array}{c}
                   y_1 \\
                   y_2 \\
                 \end{array}
               \right)
,\quad {\bf\breve{y}}:=E_2\left(
                 \begin{array}{c}
                   y_1 \\
                   y_2 \\
                 \end{array}
               \right)=E_2E_1^{-1}\mathbf{y}.
\end{equation}
It is clear that $\mathbf{y},\mathbf{\breve{y}}\in \mathbf{l}(\I)$.
On the other hand, by the invertibility of  $E_1$ and $E_2$, one can get that
for any $\mathbf{y}\in \mathbf{l}(\I)$ or $\mathbf{\breve{y}}\in \mathbf{l}(\I)$, there exist uniquely two $y_i\in l(\I)$,  $i=1,2$,
such that (\ref{y y y}) holds. In fact,
\begin{equation}\label{y1 y2 y}
y_1=(I_{2n},0)E_1^{-1}\mathbf{y}=(I_{2n},0)E_2^{-1}\mathbf{\breve{y}},\quad y_2=(0, I_{2n})E_1^{-1}\mathbf{y}=(0, I_{2n})E_2^{-1}\mathbf{\breve{y}}.
\end{equation}

Consider the following system
\begin{align}\label{HS3}
\mathbf{J}\De \mathbf{y}(t)-\mathbf{P}(t)\mathbf{R}(\mathbf{y})(t)
= \mathbf{W}(t)\mathbf{R}(\mathbf{g})(t),\quad t \in \I,
\end{align}
where  $\mathbf{W}={\rm diag}\{W_1, W_1, W_2, W_2\}$ is a diagonal matrix,
\begin{align*}
&\mathbf{J}=\left(
     \begin{array}{cccc}
       0 & 0 &-I_n & 0\\
       0 & 0 & 0& -I_n \\
       I_n &0 & 0 & 0 \\
      0& I_n & 0 & 0
     \end{array}
   \right),\quad    \mathbf{P}=\left(
     \begin{array}{cccc}
       0 & -C & D & 0 \\
      -C^*& 0 & 0& A^* \\
       D^* & 0 & 0 & B^* \\
       0&  A & B & 0
     \end{array}
   \right).
\end{align*}
$A, B, C, D$ are $n\times n$ matrices defined as (\ref{P=A,B,C,D}),
and the partial right shift operator defined on $\mathbf{l}(\I)$ is
\begin{equation*}
\mathbf{R}(\mathbf{y})(t)=\left(
                            \begin{array}{c}
                             \mathbf{u}(t+1) \\
                              \mathbf{v}(t) \\
                            \end{array}
                          \right),\quad \forall\;\mathbf{y}(t)=\left(
                            \begin{array}{c}
                             \mathbf{u}(t) \\
                              \mathbf{v}(t) \\
                            \end{array}
                          \right)\in \mathbf{l}(\I),
\end{equation*}
with  $\mathbf{u}(t)$, $\mathbf{v}(t)\in \C^{2n}$.

It is clear that $\mathbf{J}$ is the $4n\times 4n$ canonical sympletic matrix,
$\mathbf{P}$ and $\mathbf{W}$ are Hermitian and $\mathbf{W}\ge 0$.
Hence, (\ref{HS3}) is a formally selfadjoint discrete linear Hamilton system.

For any $\mathbf{x},\mathbf{y}\in \mathbf{l}(\I)$, it can be verified directly that
\begin{eqnarray}\label{yJx}
\mathbf{y}^*(t)\mathbf{J}\mathbf{x}(t)=(y_2^*Jx_1)(t)+(y_1^*Jx_2)(t), \quad t\in \I.
\end{eqnarray}
where $y_i$ and $x_i$, $i=1,2$, are determined by (\ref{y1 y2 y}).

The natural difference operator corresponding to system $(\ref{HS3})$ is
\begin{align*}
(\mathbf{L }\mathbf{y})(t)&:=\mathbf{J}\Delta \mathbf{y}(t)-\mathbf{P}(t)\mathbf{R}(\mathbf{y})(t),
\quad \mathbf{y}\in \mathbf{l}(\I).
\end{align*}
Further, similar as that of Section 2, we  define
\begin{align*}
\mathbf{l}^2_{\mathbf{W}}(\I)&:=\left\{\mathbf{y}\in \mathbf{l}(\I):\;\sum_{t\in
\I}\mathbf{R}(\mathbf{y})^*(t)\mathbf{W}(t)\mathbf{R}(\mathbf{y})(t)<+\infty\right\},\\
\langle \mathbf{y},\mathbf{z}\rangle_{\mathbf{W}}&:=\sum_{t\in \I}\mathbf{R}^*(\mathbf{z})(t)\mathbf{W}(t)\mathbf{R}(\mathbf{y})(t),\quad
\|\mathbf{y}\|_{\mathbf{W}}:=(\langle\mathbf{y}, \mathbf{y}\rangle_{\mathbf{W}})^{1/2},\quad \mathbf{y},\mathbf{z}\in \mathbf{l}^2_{\mathbf{W}}(\I),\\
\mathbf{L}^2_{\mathbf{W}}(\I)&:=\mathbf{l}^2_{\mathbf{W}}(\I)/\{\mathbf{y}\in \mathbf{l}^2_{\mathbf{W}}(\I):\; \|\mathbf{y}\|_{\mathbf{W}}=0\},\\
\mathbf{l}^2_{\mathbf{W},0}(\I)&:=\{\mathbf{y}\in \mathbf{l}^2_{\mathbf{W}}(\I):\; \exists\; s,k\in \I \;
{\rm s.\; t.}\; \mathbf{y}(t)=0
\;{\rm for}\; t\leq s\; {\rm and}\; t\ge k+1 \}.
\end{align*}
 $\mathbf{L}^2_{\mathbf{W}}(\I)$ is a Hilbert space with $\langle \cdot,\cdot\rangle_{\mathbf{W}}$.

The following results can be easily verified  and we omit the proof.

\begin{lemma}\label{equivalent 1}
Let $\mathbf{y}$, $\mathbf{\breve{y}}$, and $y_i$, $i=1,2$, have the relationship given by (\ref{y y y}) and (\ref{y1 y2 y}).
\begin{enumerate}

\item  $y_i,g_i$  satisfy  $(\Ls_iy_i)(t)= W(t)R(g_i)(t)$ on $\I$, $i=1,2$,
if and only if $\mathbf{y}$ and $\mathbf{\breve{g}}$ satisfy (\ref{HS3}).

\item $\mathbf{y}\in \mathbf{l}^2_{\mathbf{W}}(\I)$
if and only if  $y_i\in l^2_{W}(\I)$, $i=1,2$,  and
\begin{equation*}
  \|\mathbf{y}\|_{\mathbf{W}}^2=\|\breve{\mathbf{y}}\|_{\mathbf{W}}^2= \|y_1\|^2+\|y_2\|^2.
\end{equation*}

\item If  $(A_1)$ holds, then the solution of
\begin{align}\label{HS4}
    \mathbf{J}\Delta \mathbf{y}(t)-\mathbf{P}(t)R(\mathbf{y})(t)=\ld \mathbf{W}(t)R(\mathbf{y})(t),\quad t \in \I.
\end{align}
with any initial value exists uniquely.
\item  If  $(A_2)$ holds, then the definiteness condition for $(\ref{HS3})$ holds, that is,
if $(\mathbf{L }\mathbf{y})(t)=\mathbf{W}(t)(\mathbf{g})(t)$ and $\mathbf{y}(k)\ne 0$ for some $k\in \I$,
then it follows that
\begin{align*}
&\sum_{t\in \I_0}R(\mathbf{y})^*(t)\mathbf{W}(t)R(\mathbf{y})(t)\\
=&\sum_{t\in \I_0}R(y_1)^*(t)W(t)R(y_1)(t)+\sum_{t\in \I_0}R(y_2)^*(t)W(t)R(y_2)(t)>0.
\end{align*}
\end{enumerate}
\end{lemma}

The maximal, pre-minimal, and minimal subspaces corresponding to (\ref{HS3}),
which are denoted by $\mathbf{H}$, $\mathbf{H}_{00}$ and $\mathbf{H}_0$, respectively, are defined as follows:
\begin{eqnarray*}
&&\mathbf{H}:=\{([\mathbf{y}],[\mathbf{g}])\in (\mathbf{L}^2_{\mathbf{W}}(\I))^2:\;{\rm \exists\;}
\mathbf{y}\in [\mathbf{y}]\;{\rm s.t.}
\; (\mathbf{L}\mathbf{y})(t)=\mathbf{W}(t)R(\mathbf{g})(t),\; t\in\I\}, \\
&&\mathbf{H}_{00}:=\{([\mathbf{y}],[\mathbf{g}])\in \mathbf{H}: \;{\rm \exists\;}
\mathbf{y}\in [\mathbf{y}]\;{\rm s.t.}\;  \mathbf{y}\in \mathbf{l}^2_{\mathbf{W},0}(\I)\;{\rm and}\; (\mathbf{L} \mathbf{y})(t)=\mathbf{W}(t)R(\mathbf{g})(t),\; t\in \I\},\nonumber\\
&&\mathbf{H}_0:=\overline{\mathbf{H}}_{00}.
\end{eqnarray*}

Since $\mathbf{P}^*(t)=\mathbf{P}(t)$ on $\I$, $\mathbf{H}_0$ is  Hermitian, and $\mathbf{H}_0^*=\mathbf{H}$ by \cite[Thoerem 3.1]{Ren2011}.
In addition, based on (4) of Lemma \ref{equivalent 1} and Remark \ref{remark}, we can write
$(\mathbf{y},[\mathbf{g}])\in \mathbf{H}$ in the sequel for simpleness.

Let $T_i$ be any extension of $H_{i,0}$, with $H_{i,0}\subseteq T_i \subseteq H_i$, $i=1,2$.
Then
\begin{eqnarray}\label{T}
\mathbf{T}:=\{(\mathbf{y},[\mathbf{\breve{g}}])\in \mathbf{H}: \; \mathbf{y},\mathbf{\breve{g}}\; {\rm are \; defined \; by \; (\ref{y y y})\; for \; all} \; (y_i,[g_i])\in T_i\}.
\end{eqnarray}
is said to be   generated by $\{T_1,T_2\}$.
On the contrary, let $\mathbf{T}$ be any extension of $\mathbf{H}_0$  with $\mathbf{H}_0\subseteq \mathbf{T}\subseteq \mathbf{H}$.
Denote
\begin{eqnarray}\label{T1T2}
&&T_1:=\{(y_1,[g_1])\in H_1: \; y_1, g_1\; {\rm are \; defined \; by \; (\ref{y1 y2 y})\; for \; all} \;
(\mathbf{y},[\mathbf{\breve{g}}])\in \mathbf{T}\}, \nonumber\\
&&T_2:=\{(y_2,[g_2])\in H_2: \;  y_2, g_2\; {\rm are \; defined \; by \; (\ref{y1 y2 y})\; for \; all} \;
(\mathbf{y},[\mathbf{\breve{g}}])\in \mathbf{T}\}.
\end{eqnarray}
$\{T_1,T_2\}$ is said to be  generated by $\mathbf{T}$.

The following results is obtained.

\begin{lemma}\label{H and H_i}
Assume that $(A_1)$ and $(A_2)$ hold. $\mathbf{T}$, and $T_i$, $i=1,2$, are related by (\ref{T}) and (\ref{T1T2}).
\begin{enumerate}
  \item  $\mathbf{T}$ is an extension of $\mathbf{H}_0$,  with $\mathbf{H}_0\subseteq \mathbf{T}\subseteq \mathbf{H}$
if and only if $T_i$ is an extension of $H_{i,0}$ with  $H_{i,0}\subseteq T_i\subseteq H_i$, $i=1,2$.
\item  $\mathbf{T}=\mathbf{H}_{0}$ if and only if $T_i=H_{i,0}$, $i=1,2$.

\item $\mathbf{T}=\mathbf{H}$ if and only if $T_i=H_i$, $i=1,2$.
\item $\mathbf{T}^*$ is generated by $\{T_2^*,T_1^*\}$.
\end{enumerate}
\end{lemma}

\begin{proof}
(1) It follows from (\ref{y y y}) and (\ref{y1 y2 y}) that
\begin{equation*}
  \mathcal{D}(\mathbf{T})=E_1( \mathcal{D}(T_1)\times  \mathcal{D}(T_2)),\quad
  \mathcal{R}(\mathbf{T})=E_2( \mathcal{R}(T_1)\times  \mathcal{R}(T_2)).
\end{equation*}
The assertion (1)  is evident by (1) and (2) of Lemma \ref{equivalent 1}.

(2)
Assume that $T_i= H_{i,0}$ and $(y_i,[g_i])\in T_{i}$, $i=1,2$.
Then $(y_i,[g_i])\in H_{i}$ and consequently
$(\mathbf{y},[\mathbf{\breve{g}}])\in \mathbf{H}$.
For any $(\mathbf{x},[\mathbf{\breve{f}}])\in \mathbf{H}$, it follows that  $(x_i,[f_i])\in H_{i}$, $i=1,2$.
Then it follows from (\ref{yJx}) and  Theorem \ref{char H0} that
\begin{eqnarray}\label{zJy}
(\mathbf{y}^*\mathbf{J}\mathbf{x})(t)=(y_2^*Jx_1)(t)+(y_1^*Jx_2)(t)=0,\quad t=a,b+1.
\end{eqnarray}
This yields that $\mathbf{T} =\mathbf{H}_{0}$  by \cite[Theorem 3.2]{Ren2014Laa}.

On the other hand, let $\mathbf{T}=\mathbf{H}_{0}$ and $(\mathbf{y},[\mathbf{\breve{g}}])\in \mathbf{T}$. Then
(\ref{zJy}) holds for all $(\mathbf{x},[\breve{\mathbf{f}}])\in \mathbf{H}$ by $\mathbf{H}_0^*=\mathbf{H}$.
In special,  by taking some $(x_1,[f_1])\in H_{1,0}$ and any $(x_2,[f_2])\in H_{2}$,
one get  that  $(y_2^*Jx_1)(a)=0$ and  $(y_2^*Jx_1)(b+1)=0$.
Inserting these into  (\ref{zJy}) separately, one get that $(y_1^*Jx_2)(a)=0$ and $(y_1^*Jx_2)(b+1)=0$.
This yields that $T_1=H_2^*= H_{1,0}$. With a similar argument, one can show $T_2= H_{2,0}$.

(3) It can be proved with a similar argument as (2). We omit the details.

(4) By $\mathbf{S}$ denote the linear relation generated by $\{T_2^*,T_1^*\}$.
For any $(y_i,[g_i])\in T_i$ and $(x_i,[f_i])\in T_i^*$, $1,2$,
set
\begin{equation*}
  \mathbf{y}=E_1\left(
         \begin{array}{c}
           y_1 \\
           y_2 \\
         \end{array}
       \right),\quad
  \mathbf{\breve{g}}=E_2\left(
         \begin{array}{c}
           g_1 \\
           g_2 \\
         \end{array}
       \right),\quad
       \mathbf{x}=E_1\left(
         \begin{array}{c}
           x_2 \\
           x_1 \\
         \end{array}
       \right),\quad
       \mathbf{\breve{f}}=E_2\left(
         \begin{array}{c}
           f_2 \\
           f_1 \\
         \end{array}
       \right).
\end{equation*}
Then $(\mathbf{y},[\mathbf{\breve{g}}])\in \mathbf{T}$, $(\mathbf{x},[\mathbf{\breve{f}}])\in \mathbf{S}$,
and it follows that
\begin{align*}
\langle \mathbf{y},\mathbf{\breve{f}}\rangle_{\mathbf{W}}-\langle \mathbf{\breve{g}},\mathbf{x}\rangle_{\mathbf{W}} =&
\sum_{t\in\I}\mathbf{R}^*(\mathbf{\breve{f}})(t)\mathbf{W}(t)\mathbf{R}(\mathbf{y})(t)-\sum_{t\in\I}\mathbf{R}^*(\mathbf{x})(t)\mathbf{W}(t)\mathbf{R}(\mathbf{\breve{g}})(t)\\
=&\sum_{t\in\I}[R^*(f_1)(t)W(t)R(y_1)(t)-R^*(x_1)(t)W(t)R(g_1)(t)]\\
&+\sum_{t\in\I}[R^*(f_2)(t)W(t)R(y_2)(t)-R^*(x_2)(t)W(t)R(g_2)(t)]\\
=&\langle y_1,f_1 \rangle-\langle x_1,g_1\rangle+\langle  y_2,f_2 \rangle-\langle x_2,g_2\rangle=0
\end{align*}
This yields that $\mathbf{S}\subset \mathbf{T}^*$.

On the contrary, let $(\mathbf{x},[\mathbf{\breve{f}}])\in \mathbf{T}^*$ with
$\mathbf{x}=E_1(x_1^T, x_2^T)^T$ and ${\bf\breve{f}}=E_2(f_1^T, f_2^T)^T$.
Then for any  $(\mathbf{y},[\mathbf{\breve{g}}])\in \mathbf{T}$
with $\mathbf{y}=E_1(y_1^T, y_2^T)^T$ and ${\bf\breve{g}}=E_2(g_1^T, g_2^T)^T$,
it follows that
\begin{align}\label{inner y f}
\langle \mathbf{y},\mathbf{\breve{f}}\rangle_{\mathbf{W}}-\langle \mathbf{\breve{g}},\mathbf{x}\rangle_{\mathbf{W}}
=&\langle y_1,f_2 \rangle-\langle x_2,g_1\rangle+\langle  y_2,f_1 \rangle-\langle x_1,g_2\rangle=0.
\end{align}
In special, by taking $y_2(t)=g_2(t)=0$ on $\I$, one get by (\ref{inner y f}) that $\langle y_1,f_1 \rangle-\langle x_1,g_1\rangle=0$.
This yields that $(x_2,[f_2])\in T_1^*$. Similarly, one can show that $(x_1,[f_1])\in T_2^*$.
Therefore,  $\mathbf{T}^*\subset \mathbf{S}$, and consequently $\mathbf{T}^*= \mathbf{S}$.
The whole proof is complete.
\end{proof}

Since $H_{i,0}\subset H_i$, $i=1,2$, it is evident that
\begin{align*}
(H_1)_s|_{\mathcal{D}(H_{1,0})}=(H_{1,0})_{s},\quad (H_2)_s|_{\mathcal{D}(H_{2,0})}=(H_{2,0})_{s}.
\end{align*}
In addition,  the following result is obtained.

\begin{lemma} \label{A_3}
Assume that $(A_1)$ and $(A_2)$ hold. Then
\begin{equation*}
H_1(0)\cap \mathcal{N}(H_2)=\{0\},\quad H_2(0)\cap \mathcal{N}(H_1)=\{0\}.
\end{equation*}
\end{lemma}

\begin{proof}
We only show $H_1(0)\cap \mathcal{N}(H_2)=\{0\}$ and the other can be proved with a similar argument.
Let $[g]\in H_1(0)\cap \mathcal{N}(H_2)$. Then there exist $(y,[g])\in H_1$ with $[y]=0$ and $(g,0)\in H_2$.
That is,
\begin{align*}
  (\Ls_1y)(t)=W(t)R(g)(t), \quad  (\Ls_2g)(t)=W(t)R(y)(t)=0,  \quad t\in \I.
\end{align*}
If $g(k)= 0$ for some $k\in \I$, then it follows from $(A_1)$ that $g(t)=0$ for all $t\in \I$.
If $g(k)\neq 0$ for some $k\in \I$, then it follows from $(A_2)$  that
\begin{equation*}
  \sum_{t\in \I_0}R(g)^*(t)W(t)R(g)(t)>0.
\end{equation*}
This yields that $W(t)R(g)(t)\ne 0$ for some $t\in \I_0$, which together with $(\Ls_1y)(t)=W(t)R(g)(t)$ on $\I$,
implies that there exists some $k\in \I$ such that $y(k)\ne 0$.
Again by $(A_2)$ one get that
\begin{equation*}
  \sum_{t\in \I_0}R(y)^*(t)W(t)R(y)(t)>0.
\end{equation*}
This contradict with the $[y]=0$. Hence,  $g(t)= 0$ for all $t\in \I$ and consequently $H_1(0)\cap \mathcal{N}(H_2)=\{0\}$ holds.
The proof is complete.
\end{proof}

Based on the above the discussion, we can get the following result.

\begin{theorem} \label{q self 1}
Assume that $(A_1)$ and $(A_2)$ hold.
\begin{enumerate}
  \item If $T$ is a quasi self-adjoint extension of  $\{H_{1,0},H_{2,0}\}$,
then  $T^*$ is a quasi self-adjoint extension of  $\{H_{2,0},H_{1,0}\}$,
and  $\mathbf{T}$ generated by $\{T, T^*\}$ is a  self-adjoint extension of $\mathbf{H}_0$.
\item If $\mathbf{T}$  is a  self-adjoint extension of $\mathbf{H}_0$, then $\mathbf{T}$ is generated by $\{T, T^*\}$,
  where $T$ is some  quasi self-adjoint extension of  $\{H_{1,0},H_{2,0}\}$.
\end{enumerate}
\end{theorem}

\begin{proof}
(1) Let $T$ be a quasi-selfadjoint extension of $\{H_{1,0},H_{2,0}\}$. It is clear that
$T^*\in {\rm Ext}\;\{H_{2,0},H_{1,0}\}$.
It follows from (3) of Lemma \ref{dual pair} that
\begin{align*}
\frac{1}{2}\dim(\mathcal{D}(H_2)/\mathcal{D}(H_{2,0}))&=\dim(\mathcal{D}(H_1)/\mathcal{D}(T))=\dim(\mathcal{D}(H_2)/\mathcal{D}(T^*))\nonumber\\
&=\dim(\mathcal{D}(T)/\mathcal{D}(H_{1,0}))=\dim(\mathcal{D}(T^*)/\mathcal{D}(H_{2,0})).
\end{align*}
Therefore,   $T^*$ is a quasi-selfadjoint extension of $\{H_{2,0},H_{1,0}\}$ by (4) of Lemma \ref{dual pair}.

Let $\mathbf{T}$ be generated by $\{T, T^*\}$.
Then $\mathbf{T}^*$ be generated by $\{T, T^*\}$ by (4) of Lemma \ref{H and H_i}.
Therefore,  $\mathbf{T}=\mathbf{T}^*$ and consequently  $\mathbf{T}$  is a  self-adjoint extension of $\mathbf{H}_0$.

(2) If $\mathbf{T}$  is a  self-adjoint extension of $\mathbf{H}_0$, then $\mathbf{T}=\mathbf{T}^*$,
and consequently $\mathbf{T}$ is generated by $\{T, T^*\}$ for some $T\in {\rm Ext}\;\{H_{1,0},H_{2,0}\}$, and
\begin{align}\label{T/H_0}
\dim(\mathcal{D}(\mathbf{T})/\mathcal{D}(\mathbf{H_{0}}))=\dim(\mathcal{D}(\mathbf{H})/\mathcal{D}(\mathbf{T})).
\end{align}
It can be easily verified that
\begin{align}\label{T/H_0 H/T}
&\dim(\mathcal{D}(\mathbf{T})/\mathcal{D}(\mathbf{H_{0}}))=\dim(\mathcal{D}(T)/\mathcal{D}(H_{1,0}))+\dim(\mathcal{D}(T^*)/\mathcal{D}(H_{2,0})),\nonumber\\
&\dim(\mathcal{D}(\mathbf{H})/\mathcal{D}(\mathbf{T}))=\dim(\mathcal{D}(H_{1})/\mathcal{D}(T))+\dim(\mathcal{D}(H_2)/\mathcal{D}(T^*)).
\end{align}
By using $(1)$ of Lemma \ref{dual pair}, one get that
\begin{align*}
\dim(\mathcal{D}(H_2)/\mathcal{D}(H_{2,0}))=\dim(\mathcal{D}(H_1)/\mathcal{D}(H_{1,0})).
\end{align*}
This, together with  (\ref{T/H_0}) and (\ref{T/H_0 H/T}) yields that
\begin{align*}
\dim(\mathcal{D}(H_{1})/\mathcal{D}(T))=\dim(\mathcal{D}(T)/\mathcal{D}(H_{1,0})).
\end{align*}
Therefore, $T$ is a quasi-selfadjoint extension of $\{H_{1,0},H_{2,0}\}$ by (4) of Lemma \ref{dual pair}.
The whole proof is complete.
\end{proof}

\begin{corollary} \label{q self 2}
 Assume that $(A_1)$ and $(A_2)$ hold.
Then $H_{1,0}$ is a quasi self-adjoint extension of  $\{H_{1,0},H_{2,0}\}$ if and only if $\mathbf{H}_0$ is self-adjoint.
\end{corollary}

\begin{proof}
It is derived directly by  Theorem \ref{q self 1}.
The proof is complete.
\end{proof}

\section{Quasi self-adjoint extensions when $\I=[a,\infty)$}

In this section,  we consider the case that $\I=[a,\infty)$. In this case, it follows from
\cite[Theorem 4.1]{Shi2006} that
\begin{equation*}
  2n\le d_{\pm}(\mathbf{H}_0)\le 4n.
\end{equation*}
Assume that  (\ref{HS4})  is in limit point case in this section, that is $d_{\pm}(\mathbf{H}_0)= 2n$.
Some sufficient conditions for (\ref{HS4})  to be in limit point case have been given \cite{Shi2006,Ren2011}.
The following is a sufficient and necessary condition for limit point case.

\begin{lemma}\label{lpc pro} \cite[Theorem 6.15]{Shi2006}.
Let $\I=[a,\infty)$, and assume that $(A_1)$ and $(A_{2})$ hold.
Then (\ref{HS4})  is  in  limit point case if and only if
\begin{equation*}
  \lim_{t\to \infty}\mathbf{y}_2^*(t)\mathbf{J}\mathbf{y}_1(t)=0,\quad \forall\; (\mathbf{y}_i,[\mathbf{\breve{g}}_i])\in \mathbf{H},\; i=1,2.
\end{equation*}
\end{lemma}

Further, one can get the following result.

\begin{lemma}\label{y2Jy1infty=0}
Let $\I=[a,\infty)$ and assume that $(A_1)$ and $(A_2)$ hold.
If (\ref{HS4}) is in limit point case, then
$(y_2^*J y_1)(\infty)=0$ and
\begin{equation}\label{inner g y}
 \langle g_1,y_2\rangle-\langle y_1,g_2\rangle=-(y_2^*Jy_1)(a)
\end{equation}
hold for any $(y_i,[g_i])\in H_i$, $i=1,2$.
\end{lemma}

\begin{proof}
For any $(y_i,[g_i])\in H_{i}$, $i=1,2$,  it follows from Lemma \ref{H1H2} that
\begin{equation*}
 \langle g_1,y_2\rangle-\langle y_1,g_2\rangle=(y_2^*Jy_1)(\infty)-(y_2^*Jy_1)(a)
\end{equation*}
Set
\begin{equation*}
\mathbf{y}=E_1\left(
                \begin{array}{c}
                  y_1 \\
                  0 \\
                \end{array}
              \right),\quad
\mathbf{\breve{g}}=E_2\left(
                \begin{array}{c}
                  g_1 \\
                 0 \\
                \end{array}
              \right),\quad
\mathbf{z}=E_1\left(
                \begin{array}{c}
                  0 \\
                  y_2 \\
                \end{array}
              \right),\quad
\mathbf{\breve{h}}=E_2\left(
                \begin{array}{c}
                  0 \\
                  g_2 \\
                \end{array}
              \right).
\end{equation*}
Then $(\mathbf{y},[\mathbf{\breve{g}}])\in \mathbf{H}$,  $(\mathbf{z},[\mathbf{\breve{h}}])\in \mathbf{H}$.
It can be easily verified that
\begin{equation}\label{zJy1}
(\mathbf{z}^*\mathbf{J}\mathbf{y})(t)=(y_2^*J y_1)(t).
\end{equation}
Further, since (\ref{HS4}) is in  limit point case,   it follows from  Lemma \ref{lpc pro} that
\begin{equation*}
 (\mathbf{z}^*\mathbf{J}\mathbf{y})(\infty)=0.
\end{equation*}
This, together with (\ref{zJy1}), yields that $(y_2^*J y_1)(\infty)=0$  and  consequently (\ref{inner g y}) holds.
\end{proof}

\begin{theorem}\label{char H0 lp}
Let $\I=[a,\infty)$ and assume  that  $(A_1)$ and $(A_{2})$ hold.
If (\ref{HS4}) is in limit point case,  then for $i=1,2$,
\begin{align*}
H_{i,0}=\{(y_i,[g_i])\in H_i: y(a)=0\}.
\end{align*}
\end{theorem}

\begin{proof}
It is a direct consequence of  Theorem \ref{char H0} by   Lemma \ref{y2Jy1infty=0}.
\end{proof}

For any subspace $Q\subset \C^{2n}$, we denote
\begin{equation*}
  Q^*:=\C^{2n}\ominus (JQ),
\end{equation*}
and define
\begin{align*}
  H_{i,00}(Q) &:=\{(y,[g])\in H_{i}:\; y(a)\in Q\; {\rm and}\; {\rm supp}\; y \;{\rm is \; compact \; in\; \I}\},  \\
   H_{i,0}(Q) &:= \overline{H_{i,00}}(Q),\quad    H_{i}(Q) := \{(y,[g])\in H_{i}:\; y(a)\in Q\}.
\end{align*}
When $Q=\{0\}$, it is clear that
\begin{align*}
    H_{i}(Q)=H_{i,0} = \{(y,[g])\in H_{i}:\; y(a)=0\}.
\end{align*}

\begin{lemma}
Let $\I=[a,\infty)$ and assume  that  $(A_1)$ and $(A_{2})$ hold.
For any subspace $Q\subset \C^{2n}$,
\begin{equation}\label{Hi0Q}
     H_{i,0}(Q)= H_{i}(Q),\quad (H_{i,0}(Q))^*=(H_{i}(Q))^*= H_{3-i}(Q^*), \quad i=1,2.
\end{equation}
\end{lemma}

\begin{proof} Let  $Q\subset \C^{2n}$ be any subspace.
For $i=1,2$, it is clear
\begin{equation*}
H_{i,00}\subset H_{i,00}(Q)\subset H_{i}(Q)\subset H_{i}.
\end{equation*}
Since $(H_{i,00})^*= (H_{i,0})^*= H_{3-i}$ by Theorem \ref{T0*=TM+},
one can get that
$ (H_{i,00}(Q))^*\subset  H_{3-i}$ and $ (H_{i}(Q))^*\subset  H_{3-i}$.
This, together with (\ref{inner g y}) yields that
\begin{align}\label{Hi00Q}
 (H_{i,00}(Q))^* &=\{(y,[g])\in H_{3-i}:\; (y^*Jx)(a)=0,\;\forall \; (x,[f])\in H_{i,00}(Q)\},  \\
 \label{HiQ}
  (H_{i}(Q))^* &=\{(y,[g])\in H_{3-i}:\; (y^*Jx)(a)=0,\;\forall \; (x,[f])\in H_{i}(Q)\}.
\end{align}

It is evident that
\begin{align*}
  \{y(a):\; (y,[g])\in H_{i,00}(Q)\}\subset  \{y(a):\; (y,[g])\in H_{i}(Q)\}\subset Q.
\end{align*}
In addition, it follows from Corollary \ref{coro patch} that
\begin{align*}
 Q\subset \{y(a):\; (y,[g])\in H_{i,00}(Q)\}.
\end{align*}
So,
\begin{align*}
  Q=\{y(a):\; (y,[g])\in H_{i,00}(Q)\}=  \{y(a):\; (y,[g])\in H_{i}(Q)\}.
\end{align*}
Therefore, the second equality of  (\ref{Hi0Q}) is implies by (\ref{Hi00Q}) and (\ref{HiQ}).
Since both $H_{i,0}(Q)$ and $ H_{i}(Q)$ are closed,
the first equality of (\ref{Hi0Q}) can by implied by the  second one.
\end{proof}

Further, the following results are obtained.

\begin{theorem}\label{Q and H_1}
Let $\I=[a,\infty)$ and assume  that  $(A_1)$ and $(A_{2})$ hold,  and (\ref{HS4}) is in limit point case. Then
\begin{enumerate}
  \item $T\in Ext\{H_{1,0},H_{2,0}\}$ if and only if there exists $Q\subset \C^{2n}$ such that
\begin{equation*}
     T= \{(y,[g])\in H_{1}:\; y(a)\in Q\}
  \end{equation*}
  \item $T$ is a quasi self-adjoint extensions of $\{H_{1,0},H_{2,0}\}$ if and only  if
  there exists $Q\subset \C^{2n}$ with $\dim Q=n$ such that
\begin{equation*}
     T= \{(y,[g])\in H_{1}:\; y(a)\in Q\}.
  \end{equation*}
\end{enumerate}
\end{theorem}

\begin{proof}
(1) We first show for any subspace $Q\subset \C^{2n}$, $H_{1}(Q)\in Ext\{H_{1,0},H_{2,0}\}$.
It is clear that
\begin{equation*}
H_{1,00}\subset H_{1,00}(Q)\subset H_{1}(Q)\subset H_{1}.
\end{equation*}
Since $ H_{1}(Q)$ and $ H_{1}$ are closed, and $(H_{1,00})^*= (H_{1,0})^*= H_{2}$ by Theorem \ref{T0*=TM+},
one can get that
\begin{equation*}
H_{1,0}=\overline{H_{1,00}}\subset H_{1}(Q)\subset H_{1}=H_{2,0}^*.
\end{equation*}
Therefore, $H_{1}(Q)\in Ext\{H_{1,0},H_{2,0}\}$ holds.

On the other hand, for any $T\in Ext\{H_{1,0},H_{2,0}\}$, set
\begin{align*}
  Q=\{y(a):\; (y,[g])\in T\}.
\end{align*}
Since $H_{1,0}\subset T\subset H_{2,0}^*=H_{1}$, it follows that $T\subset H_{1}(Q)$.

For any $(y,[g])\in H_{1}(Q)$, there exists $(x,[f])\in T$ such that $y(a)=x(a)$.
Set $(z,[h])=(y,[g])-(x,[f])$. Then $(z,[h])\in H_{1}$ and $z(a)=0$.
It follows from Theorem \ref{char H0 lp} that $(z,[h])\in H_{1,0}\subset T$.
Therefore, $(y,[g])=(z,[h])+(x,[f])\in T$. By the arbitrary of  $(y,[g])\in H_{1}(Q)$
one get that $H_{1}(Q)\subset T$. Hence, the first assertion holds.

(2) If $ H_{1}(Q)\in Ext\{H_{1,0},H_{2,0}\}$ is a quasi-selfadjoint extension, then
\begin{equation}\label{dim}
\dim(\mathcal{D}(H_1(Q)/\mathcal{D}(H_{1,0}))=\dim(\mathcal{D}(H_1(Q))^*)/\mathcal{D}(H_{2,0}))=\dim(\mathcal{D}(H_2(Q^*))/\mathcal{D}(H_{2,0})).
\end{equation}
where the second equality of (\ref{Hi0Q}) is used.
It is clear
\begin{equation}\label{dim1}
\dim(\mathcal{D}(H_1(Q)/\mathcal{D}(H_{1,0}))=\dim Q, \quad
\dim(\mathcal{D}(H_2(Q^*))/\mathcal{D}(H_{2,0}))=\dim Q^*.
\end{equation}
This, together with the fact $\dim Q+\dim Q^*=2n$, yields that $\dim Q=n$.

On the contrary, if $\dim Q=n$, then $\dim Q=\dim Q^*=n$, and consequently
$(\ref{dim1})$ holds. Again by (\ref{Hi0Q}), one get $H_1(Q)=H_2(Q^*)$.
So, $(\ref{dim})$ holds, which implies that $ H_{1}(Q)\in Ext\{H_{1,0},H_{2,0}\}$ is a quasi-selfadjoint extension.
The proof is complete.
\end{proof}

\section*{Acknowledgement}
This work is supported by the NSF of Shandong Province, P.R. China [grant numbers  ZR2020MA012].

\bibliographystyle{amsplain}

\end{document}